\documentclass[ijoo,nonblindrev]{informs-ijoo}

\OneAndAHalfSpacedXI

\usepackage{natbib}
 \bibpunct[, ]{(}{)}{,}{a}{}{,}%
 \def\bibfont{\small}%
\TheoremsNumberedThrough     
\ECRepeatTheorems

\EquationsNumberedThrough    

\usepackage{graphicx}
\usepackage{algorithm,algpseudocode}
\usepackage{mathrsfs}
\usepackage{amsmath}  
\usepackage{mathtools}
\usepackage{amssymb}
\usepackage{amsfonts} 
\usepackage{subcaption}
\usepackage[export]{adjustbox}
\usepackage{bbm}
\usepackage{cite}
\usepackage{url}
\usepackage{color}
\allowdisplaybreaks[1]

\newtheorem{observation}{Observation}
\newcommand{\norm}[1]{\left\lVert#1\right\rVert}
\DeclarePairedDelimiterX\Set[2]{\lbrace}{\rbrace}{ #1 \,\delimsize| \,\mathopen{} #2 }

\newcommand{\mtc}{\mathcal}
\newcommand{\mb}{\mathbb}

\newcommand{\cP}{\mathcal{P}}

\newcommand{\cS}{\mathcal{S}}
\newcommand{\cG}{\mathcal{G}}
\newcommand{\F}{\mathcal{F}}

\newcommand{\ind}[1]{\mathbbm{1}\{#1\}}
 
\newcommand{\wh}[1]{\widehat{#1}}

\newcommand{\bs}[1]{\boldsymbol{#1}}
\newcommand{\E}{\mathbb{E}}
\newcommand{\var}{\mathbb{V}\textrm{ar}}

\newcommand{\diag}{\textrm{diag}}
\renewcommand{\exp}{\textrm{exp}}

\newcommand{\beq}{\begin{equation}}
\newcommand{\eeq}{\end{equation}}
\newcommand{\ba}{\begin{aligned}}
\newcommand{\ea}{\end{aligned}}
\newcommand{\bdm}{\begin{displaymath}}
\newcommand{\edm}{\end{displaymath}}

\begin{document}



\RUNTITLE{Service Center Location Problem with Decision Dependent Demand}

\TITLE{A Distributionally-Robust Service Center Location Problem with  
Decision Dependent Demand Induced from a Maximum Attraction Principle}

\ARTICLEAUTHORS{%
\AUTHOR{Fengqiao Luo}
\AFF{Department of Industrial Engineering and Management Science, Northwestern University, Evanston, Illinois 60208, \EMAIL{fengqiaoluo2014@u.northwestern.edu}} 
} 

\ABSTRACT{%
This paper establishes and analyzes a service center location model with a simple but novel decision-dependent
demand induced from a maximum attraction principle. The model formulations are investigated
in the distributionally-robust optimization framework for the capacitated and uncapacitated cases. A statistical model that is
based on the maximum attraction principle for estimating
customer demand and utility gain from service is established and analyzed.
The numerical experiments show that the model admits high computational efficiency
in solving mid- and large-size instances.
}%


\maketitle

\section{Introduction}
The facility location problem is one 
of the most fundamental problems investigated in operations research.
In this problem, a decision maker needs to decide locations 
of a limited number of  facilities (factories, retail centers, power plants, service centers, etc.), 
and determine coverage of demand from different sites by the located facilities. 
The objective is to minimize the facility setup cost and the cost of production and delivery. 
This problem provides a basic framework to formulate related problems in resource allocation, 
supply chain management and logistics, etc.
The facility location models in the stochastic and robust (distributionally-robust) optimization framework
have also received sufficient investigation. In the stochastic optimization framework,
the customer demand can be random parameters following a known or partially known 
probability distributions, while in the robust (distributionally-robust) optimization
framework customer demand (the probability distribution of customer demand) 
can have a certain level of uncertainty (ambiguity). Incorporating randomness
or uncertainty in customer demand estimation is a sensible and realistic model improvement,
as the decision maker has no perfect prediction on the demand in practice. 
Recently, there is a trend of research on distributionally-robust optimization problems
with decision-dependent ambiguity, in which the ambiguity set is specified 
by decision dependent parameters. Introducing decision dependency
has certain merit in situations when model parameters are naturally functions of decision variables.
For example, in a basic pricing and revenue management problem 
the demand of a product can be a decreasing function of the price to be determined.  
The notion of modeling decision-dependency in the robust (distributionally-robust) optimization
framework can be applied to investigate a variety of specific problems.

In this paper, we investigate a distributionally-robust
service center location problem with decision-dependent customer demand.
In the model, customers need to physically access a service center
in order to receive service. This problem is motivated from clinic or medical
test center (such as for COVID-19 screen test) allocation.
The objective is to maximize the total utility gain of all customers who
have decided to receive service from the opened service centers.
The decision dependence is decoupled from the ambiguity set in our model,
which is different from \citep{basciftci2020-dro-fl-endog-dmd}.
Specifically, we impose the decision dependency on the demand itself
but not on the uncertainty of demand, and the decision-dependent demand
is induced from a maximum attraction principle based on a ranking 
of opened service centers in the neighborhood. 
Our investigation shows that the decoupling approach leads to
high computational efficiency of solving the problem and it can
be more data-driven in practice. Furthermore, the novel approach
of modeling decision dependency introduce additional combinatorial
properties to the model which is of independent interest.

The contribution of this paper is summarized as follows:
\begin{enumerate}
	\item A novel approach of modeling decision dependency of the customer demand
		has been established for the service center location problem, with a possible extension
		on the modeling of decision dependency. In this approach, the decision dependency is decoupled
		from the ambiguity set which admits highly computational tractable 
		reformulations.
	\item A learning model has been established to estimate utilities and demand based on
		data from survey.
	\item The numerical experience with this model shows that mid- and large-size instances can be solved
		very efficiently due to decoupling of decision dependency from the ambiguity set.
\end{enumerate}

\subsection*{Literature review}
Facility location problems with non-deterministic demand have received 
plentiful investigation in the framework of stochastic optimization and 
robust (distributionally-robust) optimization. In the stochastic optimization
framework, customer demand are independent random parameters with known probability
distributions and the problem can be formulated as
a two-stage stochastic program in which the location vector is in the 
first-stage decision that needs to be made before realization of 
customer demand \citep{1992-dual-stoch-fac-loc,2011-fac-loc-bern-dmd}. 
In a special case when the customer demand rates and service rates of each facility 
are assumed to follow exponential probability distributions,
the problem of minimizing long term average cost can be reformulated as   
a deterministic optimization problem with corresponding Poisson rates
to characterize the demand and service. This setting has been applied
to an allocation problem of ATMs \citep{2002-alg-fac-loc-stoch-demand}.  
The probability distribution of demand can also be used to define 
chance constraints to ensure a required service level under possible
stockout and supply disruption
\citep{2011-rob-fac-loc-bio-terror-attack,gulpinar2013-rob-fac-loc,
daskin2012_fac-loc-rand-disrupt-imp-estim,2017-multi-source-suppl-chance-constr}.

In the robust (distributionally-robust) optimization framework
for facility location problems,
the information of demand is partially known to the decision maker,
and the goal is to find a robust optimal location vector that optimizes
the objective after a worst-case realization of demand information
\citep{2006-fac-loc-uncert-rev}. 
Following this direction of research,  
\citealt{2010-fac-loc-rob-opt-approach} investigated a robust multi-period facility location problem
with box and ellipsoid sets of uncertainty.
\citealt{2020-gourtani_dro-two-stage-fac-loc} 
investigated a distributionally-robust two-stage facility location problem
with an ambiguity set defined corresponding to the mean and covariance matrix
of a random parameter for expressing the demand.
Facility location problems with demand uncertainty have been 
studied with a variety of novel application background, which includes
but not limited to medication coverage and delivery under a large-scale bio-terror attack
\citep{2011-rob-fac-loc-bio-terror-attack}, medical equipment (defibrillators) location 
problem to reduce cardiopulmonary resuscitation (CPR) risk \citep{2017-rob-deploy-card-arrest-loc-uncert},
humanitarian relief logistics \citep{2012-two-echelon-stoch-fac-loc}, 
and hazardous waste transportation \citep{2014-rob-fac-loc-hazard-waste-transport}, etc.

In a lot of real world problems, the uncertainty of parameters can interplay with
the decision to be made. This behavior is well observed especially in a sequential (multi-stage) decision-making process,
in which information about system parameters are gradually revealed and the decisions made 
up until the current stage can reshape the uncertainty in future \citep{grossmann2006_stoch-prog-dec-uncert}.
It motivates the research on multi-stage stochastic optimization with decision-dependent uncertainty.
Solution strategies based on Lagrangian duality and novel branch-and-bound
methods are developed to solve this family of problems 
\citep{goel2005_lag-dual-BB-lin-stoch-prog-ddu,
gupta2011_sol-strat-multi-stoch-prog-endog-uncert,
tarhan2013_comp-methd-nonconv-multistage-MINLP-dec-dept-uncert},
and it has a rich application in oil-chemical industrial 
\citep{grossmann2004_stoch-prog-plan-gas-field-uncert-reserv,
tarhan2008_multi-sddd-oil-synthesis,goel2006_BB-alg-gas-field-reserv-uncert,
grossmann2009_stoch-prog-plan-oil-infrast-dec-uncert}.
An approximation scheme \citep{vayanos2011-dec-rule} is proposed to tackle the high complexity
in solving the multi-stage problems with decision-dependent
information discovery. Decision-dependent uncertainty
has also been considered in stochastic optimization problems 
such as resource management \citep{tsur2004_gdwater-threat-catastr-event}, 
stochastic traffic assignment modeling \citep{shao2006_reliab-stoch-traff-assign-demand-uncert}, 
and robust network design \citep{ahmed2000-phd_plan-uncert-stoch-MIP,
viswanath2004_invest-stoch-netwk-min-exp-short-path}.

Imposing decision-dependent uncertainty for robust (distributionally-robust) optimization
has received great attention in recently years 
\citep{nohadani2016_opt-dec-uncert,luo-2018-D3RO,noyan2018_dro-dec-dept-amb-set,
basciftci2020-dro-fl-endog-dmd,luo2019_DRO-FL}.
\citealt{nohadani2016_opt-dec-uncert} studied robust 
linear programs with decision dependent budget-type uncertainty
and its generalization with a polyhedral uncertainty set. 
This concept is demonstrated in a robust shortest-path problem,
where the uncertainty is resolved progressively when approaching the destination.  
\citealt{noyan2018_dro-dec-dept-amb-set} investigated a family of distributionally-robust
optimization problems with an ambiguity set defined using earth mover's distances
(including total variation distance and the Wasserstein metric) 
with decision-dependent parameters such as the nominal probability distribution 
and the radius, and focused on understanding which settings can lead to
tractable formulation. \citealt{royset2017_var-thy-opt-u-stoch-amb} provided
a variational principle analysis for optimization under stochastic ambiguity,
which gives fruitful tools for analyzing solution quality and price of robustness.
Some novel applications are studied in
\citep{spacey2012_rob-software-partition-multi-instant} for robust software
partition, and in \citep{nohadani2017-rob-opt-rad-thpy-time-dept-uncert} for radiation
therapy design. 

The work in \citep{basciftci2020-dro-fl-endog-dmd} and 
\citep{luo2019-serv-loc-dec-dept-utility} are mostly related to the work of this paper.
\citealt{basciftci2020-dro-fl-endog-dmd} investigated
a distributionally-robust facility location problem with decision-dependent
ambiguity, where the ambiguity set is defined using disjointed lower and upper bounds on
the mean and variance of the candidate probability distributions for customer demand,
and these bounds are defined as linear functions of the location vector to admit
a tractable MILP reformulation.
The model considered in \citep{luo2019-serv-loc-dec-dept-utility} assumed 
deterministic customer demand and imposed a linear decision dependency on the
moment bounds for candidate probability distributions of utility, which leads to a MISOCP reformulation.

\section{Problem Formulation}
We consider a problem of allocating a set of service centers (facilities)
in a region to meet the customer demand from the region. Similarly to the traditional
facility location problem, the region consists of customer sites and candidate
locations for the service centers. Each customer site has certain demand that
needs to be fulfilled by a service center from its neighborhood.
Especially, we consider the case that the demand of each customer site 
is a random parameter which depends on the locations of service centers. 
In most literature on decision-dependent robust (distributionally-robust) optimization,
the decision-dependency is formulated as a linear function. 
In this problem, we consider a special form of decision dependency
that is based on an assumption of maximum attraction principle,
which will be clear in Section~\ref{sec:determ-model}. 

\subsection{A deterministic model}\label{sec:determ-model}
We first formulate a deterministic model of this problem, in which
we treat the demands as deterministic parameters that depend on 
the locations of service centers and introduce the notion of maximum
attraction principle. The notations used in this model is given by the
following list:
\begin{table}[H]
\centering
\begin{tabular}{l|l}
\hline\hline
$\cS$ & the set of customer sites; \\
$\F$ & the set of candidate locations for service centers; \\
$b_j$ & the cost of opening a service center at location $j\in\F$; \\
$B$ & the budget of allocating service centers; \\
$C_j$ & the service capacity of a service center located at $j\in\F$; \\
$u_{ij}$ & the utility gain obtained by a customer from site $i\in\cS$ \\
	     &	 who gets service from the service center at $j\in\F$; \\
$y_j$ & the binary decision variable of opening a service center at location $j\in\F$; \\ 
$D_i(y)$ & the demand from customer site $i$ that depends on the decision vector $y$; \\
$x_{ij}$ & the demand flow from customer site $i\in\cS$ to the service center at $j\in\F$. \\
\hline
\end{tabular}
\end{table}
The deterministic model is formulated as follows:
\beq\label{opt:FL-determ1}
\begin{aligned}
&\max\;\sum_{i\in\cS}\sum_{j\in\F} u_{ij}x_{ij} \\
&\textrm{ s.t. }\;\sum_{j\in\F}b_jy_j\le B, \\
&\qquad\; \sum_{i\in\cS}x_{ij} \le C_jy_j \qquad\forall j\in\F, \\
&\qquad\; \sum_{j\in\F}x_{ij} \le D_i(y) \qquad \forall i\in\cS, \\
&\qquad\; y_j\in\{0,1\},\; x_{ij}\ge 0 \;\forall i\in\cS, j\in\F.
\end{aligned}
\eeq
The first three constraints represent the limit of budget, capacity and 
the amount of demand, respectively. Note that in general the demand from
site $i\in\cS$ can depend on the pattern of service center locations which
is represented by the decision vector $y$. 
In this paper, we investigate 
a specific type of location dependent demand which is simple but sensible
in depicting customers' behavior in practice. This type of location dependency of demand is
referred as the \emph{maximum attraction principle} in this paper.
We now give a detailed description of this principle. 
First, each customer site $i$ is associated with a preferable subset $\F_i$
of candidate facility locations. When there is only one service center opened 
at $j\in\F_i$, it will attract $D_{ij}$ amount
of demand from $i$. If there are multiple service centers opened at 
locations in $\F_i$, the maximum demand that can be attracted from $i$ 
is equal to the maximum $D_{ij}$ for $j\in\F_i$ that has a facility. 
This principle is formally described by the following equation:
\beq\label{eqn:max-attr}
D_i(y)=\left\{
\begin{array}{ll}
0 & \textrm{ if } \sum_{j\in\F} \ind{j\in\F_i} y_j=0 \\
D_{ij} & \textrm{ if } y=e_j \;\textrm{ for some } j\in\F_i \\
\max_{j\in\F_i(y)} D_{ij} & \textrm{ for other cases, } 
\end{array}
\right.
\eeq
where $e_j$ is the $|\F_i|$-dimensional vector with the $j$-th entry being 1
and other entries being 0, and the set $\F_i(y)$ is defined as $\F_i(y):=\Set*{j\in\F_i}{y_j=1}$. 
The maximum attraction principle matches with our intuition from practice.
This principle first assumes that locations of service centers that are not
within the preferable location set $\F_i$ are not attractive to customers
from site $i$ at all. The most natural way of establishing $\F_i$ is based on
the distance from site $i$ to the candidate locations. 
This is a reasonable assumption because customers usually will not consider visiting
a service center that is beyond a certain distance from their living place.
The maximum attraction principle further assumes that within the preferable
opened service centers $\F_i(y)$, there exists one service center $j^*$ (or multiple centers) that  
is (are) most attractive to the customers from site $i$ in the sense that 
$j^*\in\textrm{argmax}_{j\in\F_i(y)}D_{ij}$, and the presence of multiple service
centers in $\F_i(y)$ including $j^*$ attracts the same amount of demand from site $i$
as the presence of just a single service center at $j^*$. This assumption is also consistent
with our intuition from practice. If every service center is identical in the sense of
scale and service quality, the most attractive one to site $i$ is likely the one
that is most close to $i$. Furthermore, customers who are willing to visit further
service centers are also willing to visit the one that is most close to their living place.
Note that the maximum attraction principle can be extended to incorporate the impact
of the number of opened service centers in $\F_i$ on the demand. 
This extension is discussed in Section~\ref{sec:model-remark}.
We also define the subset $\cS_j$ as $\cS_j:=\Set*{i\in\cS}{j\in\F_i}$ for all $j\in\F$.
Based on the maximum attraction principle, the service center location problem 
\eqref{opt:FL-determ1} is written as the follows:
\beq\label{opt:FL-determ2}
\begin{aligned}
&\max\;\sum_{i\in\cS}\sum_{j\in\F} u_{ij}x_{ij} \\
&\textrm{ s.t. }\;\sum_{j\in\F}b_jy_j\le B, \\
&\qquad\; \sum_{i\in\cS}x_{ij} \le C_jy_j \qquad\forall j\in\F, \\
&\qquad\; \sum_{j\in\F_i}x_{ij} \le\max_{j\in\F_i(y)} D_{ij} \qquad \forall i\in\cS, \\
&\qquad\; y_j\in\{0,1\},\; x_{ij}\ge 0 \;\forall i\in\cS, j\in\F.
\end{aligned}
\eeq 
To go one step further, we can linearize the term $\max_{j\in\F_i(y)} D_{ij}$ by introducing
some continuous auxiliary variables $q_{ij}$ to form it as a convex combination
of $D_{ij}$ for $j\in\F_i$. After this transformation, we obtain the following 
equivalent formulation of \eqref{opt:FL-determ2}:
\begin{subequations}
\makeatletter
        \def\@currentlabel{DDSL}
        \makeatother
        \label{opt:DDSL}
        \renewcommand{\theequation}{DDSL.\arabic{equation}}
\begin{align}
&\max\;\sum_{i\in\cS}\sum_{j\in\F} u_{ij}x_{ij} \\
&\textrm{ s.t. }\;\sum_{j\in\F}b_jy_j\le B,  \label{constr:determ1}\\
&\qquad\; \sum_{i\in\cS}x_{ij} \le C_jy_j \qquad\forall j\in\F,  \label{constr:determ2} \\
&\qquad\; \sum_{j\in\F_i}x_{ij} \le\sum_{j\in\F_i}D_{ij}q_{ij} \qquad \forall i\in\cS,  \label{constr:determ3}\\
&\qquad\; \sum_{j\in\F_i}q_{ij}\le 1 \qquad \forall i\in\cS,  \label{constr:determ4}\\
&\qquad\; q_{ij}\le y_j \qquad \forall i\in\cS,\forall j\in\F_i,  \label{constr:determ5}\\
&\qquad\; y_j\in\{0,1\},\; x_{ij}\ge 0, \; q_{ij}\ge 0\;  \forall i\in\cS, j\in\F_i.  \label{constr:determ6}
\end{align}
\end{subequations}
Notice that in the above reformulation, the auxiliary variables 
$\Set*{q_{ij}}{j\in\F_i}$ are used select the most attractive candidate
location driven by the sense of maximizing the objective.
The constraint $q_{ij}\le y_j$ ensures that only opened service
centers in $\F_i$ are involved in the maximum attraction principle.
Since $\sum_{j\in\F_i}q_{ij}\le1$, the model will set $q_{ij^*}=1$
and $q_{ij^\prime}=0$ for all $j^\prime\in\F_i\setminus\{j^*\}$
to relax the constraint $\sum_{j\in\F_i}x_{ij} \le\sum_{j\in\F_i}D_{ij}q_{ij}$
as possible.

\subsection{Some remarks on the model}\label{sec:model-remark}
In the model setting, we implicitly assume that the customers are willing 
to corporate with the decision maker to maximize the total utility gain. 
Although this assumption is highly impractical,
the rationality of this model depends on what metric we use to measure
the system performance. If the goal is to estimate what is expectation
of total utility in practice, then the following questions should be addressed 
and the corresponding aspects should be properly modeled:
What is the service policy used by each service center (FIFO or some other policies)?
How to characterize customers' behavior and the mechanism of competing for
the limited service capacity? If their most favored service center does not have
any capacity, are they willing to accept the service from the less preferred locations
with less utility gain?
Incorporating all these factors into the model can easily make it very complicated, 
and meanwhile a large amount of customers' information are needed to drive
this approach of modeling. On the other side, if the goal is to access what is the
maximum potential utility that can be achieved by the system in the most ideal situation,
then the \eqref{opt:DDSL} model can be used to give an estimation. 
A possible modification of modeling the decision-dependent demand
is to add a perturbation term to the demand based on the number of 
opened service centers. Specifically, we can modify the constraint
$\sum_{i\in\F_i}x_{ij}\le\max_{j\in\F_i}D_{ij}$ to be
$\sum_{i\in\F_i}x_{ij}\le\max_{j\in\F_i}D_{ij}+a_i(\sum_{j\in\F_i}y_j)$,
where $a_i$ is a parameter that measures the influence of opening
one more service center on the demand. In this way of modeling,
it is assumed that the marginal increase of demand may depend
on the number of opened service centers. In reality, when potential 
customers see more chain stores are opened in the neighborhood,
they may have higher intention to try one of them. In this case,
the parameter $a_i$ can be a random parameter with certain level of ambiguity
in the distributionally-robust extension of the model.

\subsection{A  distributionally-robust two-stage  stochastic extension}
The deterministic service center location model with decision dependent
demand can be further extended to a distribtuionally-robust two-stage 
stochastic program after imposing an ambiguity set on the pairwise 
demand $D_{ij}$. This extension is motivated by the fact that 
estimation of the demand parameters $D_{ij}$ could be inaccurate. 
Therefore, we can treat $D_{ij}$ as random parameters with an unknown
joint probability distribution, and apply the distributionally-robust framework
on this service center location problem with uncertainty. We assume
a finite support of the joint probability distribution of the pairwise random demand
vector $D:=\Set*{D_{ij}}{i\in\cS,\;j\in\F_i}$. Let the finite support
be based on $|\Omega|$ samples written as 
$D^\omega=\Set*{D^{\omega}_{ij}}{i\in\cS,\;j\in\F_i}$ for all $\omega\in\Omega$.
In this case, any probability distribution of $D$ can be represented as
a $|\Omega|$-dimensional vector. We define a nominal probability distribution 
$\mu_0$ of $D$ as 
\beq
\mu_0(D=D^\omega)=\mu^{\omega}_0 \qquad \forall\omega\in\Omega.
\eeq
In the vector representation, we write it as $\mu_0=[\mu^{\omega}_0:\omega\in\Omega]$.
The ambiguity set of candidate joint probability distribution of $D$ is defined
based on the total variation distance between two probability distributions.
Specifically we consider an ambiguity set of the following form:
\beq\label{eqn:amb-set}
\cP:=\Set*{\mu\in\mb{R}^{|\Omega|}}{\norm{\mu-\mu_0}_1\le d}.
\eeq
Based on the above definition of ambiguity set,
the distributionally-robust two-stage stochastic extension of the service
center location model can be formulated as follows:  
\beq\label{opt:FL-DRO}
\begin{aligned}
&\underset{y}{\max}\;\underset{\mu\in\cP}{\min}\;\E_{D\sim\mu}[Q(y,D)] \\
&\textrm{ s.t. }\; \sum_{j\in\F} b_jy_j\le B, \\
&\qquad\; y_j\in\{0,1\}\;\forall j\in\F,
\end{aligned}
\tag{DRO-FL}
\eeq
where the recourse function $Q(y,D^\omega)$ for scenario $\omega$ 
is given by
\begin{subequations}\label{opt:second-stage}
\begin{align}
Q(y,D^\omega)=&\max\;\sum_{i\in\cS}\sum_{j\in\F_i} u_{ij}x^{\omega}_{ij} \\
&\textrm{ s.t. }\; \sum_{i\in\cS} x^{\omega}_{ij}\le C_jy_j \qquad\forall j\in\F, \label{cst:Q1} \\
&\qquad\; \sum_{j\in\F_i}x^{\omega}_{ij}\le\sum_{j\in\F_i}D^{\omega}_{ij}q^{\omega}_{ij} \qquad\forall i\in\cS, \label{cst:Q2}\\
&\qquad\; \sum_{j\in\F_i}q^{\omega}_{ij}\le 1 \qquad\forall i\in\cS, \label{cst:Q3}\\
&\qquad\; q^{\omega}_{ij}\le y_j \qquad\forall i\in\cS,\forall j\in\F_i, \label{cst:Q4}\\
&\qquad\; x^{\omega}_{ij}\ge 0,\;q^{\omega}_{ij}\ge 0\;\forall i\in\cS,\;\forall j\in\F_i.
\end{align}
\end{subequations}
Using a standard technique for two-stage stochastic programming,
we can decompose \eqref{opt:FL-DRO} into a master problem and 
scenario sub-problems. The original problem \eqref{opt:FL-DRO}
can be solved iteratively. In each iteration, we solve the current master
problem in the space of $y$ and pass the current master solution
to each scenario sub-problem. After solving each scenario sub-problem
for the fixed first-stage solution, we can generate a valid inequality for
each scenario using the optimal dual values associated with constraints
of the scenario sub-problem. Then we aggregate the valid inequalities 
from all scenario sub-problems using the worst-case measure of scenarios
to get a single cut for the master problem which is an optimality cut. 
The optimality cut is added to the master problem in the next iteration. 
Specifically, the master problem at iteration $n$ can be represented as
\beq\label{opt:master}
\ba
&\max\; \eta \\
&\textrm{ s.t. } \sum_{j\in\F} b_jy_j\le B, \\
&\qquad\;\eta \le \sum_{\omega\in\Omega}\mu^{(k)}_{\omega}\Big(r^{(k)}_{\omega}+\sum_{j\in\F}t^{(k)}_{\omega,j}y_j \Big)
\qquad \forall k\in[n-1], \\
&\qquad y\in\{0,1\}^{|\F|},\;\eta\ge 0,
\ea
\tag{Master}
\eeq
where $\mu^{(k)}$ is the iteration based worst-case probability measure on scenarios 
at iteration $k\in[n-1]$. The way of determining this worst-case probability measure
is given in \eqref{opt:LP-worst-prob}.
At iteration $n$ we solve the master problem \eqref{opt:master}
and obtain the current optimal first-stage solution $y^{(n)}$.
This $y^{(n)}$ is input into every second-stage scenario sub-problem.
The second-stage linear program is solved and let $\alpha^{(n)}_{\omega,j}\ge 0$, 
$\beta^{(n)}_{\omega,i}\ge 0$, $\gamma^{(n)}_{\omega,i}\ge 0$, 
and $\tau^{(n)}_{\omega,ij}\ge0$ 
for all $i\in\cS,\;j\in\F_i$ are optimal dual values corresponding to the constraints 
\eqref{cst:Q1}, \eqref{cst:Q2}, \eqref{cst:Q3} and \eqref{cst:Q4}, respectively.
Using the standard technique of Bender's decomposition [cite ref],
we can obtain the following valid inequality on the value function $Q(y,D^\omega)$:
\beq
Q(y,D^\omega)\le \sum_{i\in\cS}\gamma^{(n)}_{\omega,i}
+\sum_{j\in\F}\big(C_j\alpha^{(n)}_{\omega,j}+\sum_{i\in\cS_j}\tau^{(n)}_{\omega,ij}\big)y_j.
\eeq
The strong duality implies that when evaluating $Q(y,D^\omega)$ at $y^{(n)}$,
we get
\beq
Q(y^{(n)},D^\omega)= \sum_{i\in\cS}\gamma^{(n)}_{\omega,i}
+\sum_{j\in\F}\big(C_j\alpha^{(n)}_{\omega,j}+\sum_{i\in\cS_j}\tau^{(n)}_{\omega,ij}
\big)y^{(n)}_j.
\eeq
The worst-case probability measure $\mu^{(n)}$ is an optimal solution of the following linear program:
\beq\label{opt:LP-worst-prob}
\ba
&\min_{\mu}\;\sum_{\omega\in\Omega}\mu_{\omega}Q(y^{(n)},D^\omega) \\
&\textrm{ s.t } \norm{\mu-\mu_0}_1\le d, \\
&\qquad \sum_{\omega\in\Omega}\mu_{\omega}=1, \; \mu_{\omega}\ge 0\;\forall\omega\in\Omega,
\ea
\eeq
where $Q(y^{(n)},D^\omega)$ is the value function evaluated at the first stage solution $y^{(n)}$
and scenario $\omega$. Once the current worst-case probability measure is obtained,
the following inequality will be added to the first-stage master problem:
\beq\label{eqn:agg-ineq}
\eta\le\sum_{\omega\in\Omega}\mu^{(n)}_{\omega}\Big[\sum_{i\in\cS}\gamma^{(n)}_{\omega,i}
+\sum_{j\in\F}\big(C_j\alpha^{(n)}_{\omega,j}+\sum_{i\in\cS_j}\tau^{(n)}_{\omega,ij}\big)y_j \Big].
\eeq
The algorithm and convergence property for solving \eqref{opt:FL-DRO}
are given in Appendix~\ref{app:cut-plane-alg}. It can be shown that 
for a given first-stage solution $y$, the scenario sub-problem can be
solved using a greedy algorithm. 

\subsection{A single-stage reformulation of \eqref{opt:FL-DRO} for the uncapacitated case}
We consider a special case of \eqref{opt:FL-DRO} in which each candidate
service center has sufficient capacity for service. In this case, the capacity constraints
\eqref{cst:Q1} can be removed from the second-stage scenario problems.
We will show that this leads to a simplified scenario problem that admits
a closed form optimal solution and optimal objective for the second-stage problem
with a mild regularity condition on the parameters.
Then after dualizing the inner minimization problem over the probability measure on 
scenarios, \eqref{opt:FL-DRO} can be reformulated as a mixed 0-1 linear program.
The single-stage reformulation result is given by the following theorem.
\begin{definition}\label{def:consist}
The utility gain and demand are \textit{consistent} if for every $\omega\in\Omega$,
$i\in\cS$ and $\F^\prime\subseteq\F_i$ there exists a $j\in\F^\prime$ such that
$u_{ij}=\max_{j^\prime\in\F^\prime}u_{ij^\prime}$ and 
$D^{\omega}_{ij}=\max_{j^\prime\in\F^\prime}D^{\omega}_{ij^\prime}$.
\end{definition}
The consistency condition for utility gain and demand says that the location
that attracts the most demand over other locations should also correspond
to the highest utility.
\begin{theorem}
Suppose the utility gain and demand are consistent.
In the uncapacitated case, the distributionally-robust service center location
problem \eqref{opt:FL-DRO} with decision dependent demand based on 
the maximum attraction principle and the ambiguity set \eqref{eqn:amb-set}
can be reformulated as the following mixed 0-1 linear program:
\beq\label{opt:FL-DRO-uncap}
\ba
&\max\; \sum_{\omega\in\Omega}\mu^{\omega}_0(\alpha^{\omega}-\beta^{\omega})-\lambda d-\gamma \\
&\;\emph{ s.t. } \sum_{j\in\F} b_jy_j\le B, \\
&\qquad\; \sum_{i\in\cS}\sum_{j\in\F_i}u_{ij}D^{\omega}_{ij}s_{ij} +\alpha^{\omega} 
-\beta^{\omega}+\gamma \ge 0 \qquad \forall\omega\in\Omega, \\
&\qquad\;\; s_{ij}\le y_{j} \qquad \forall i\in\cS, \forall j\in\F_i, \\
&\qquad \sum_{j\in\F_i}s_{ij}\le 1 \qquad \forall i\in\cS, \\
&\qquad\; \lambda\ge0,\;\alpha^{\omega}\ge0,\;\beta^{\omega}\ge0,\;\gamma\in\mathbb{R},\;\forall \omega\in\Omega\\
&\qquad\; y_j\in\{0,1\}\;\forall j\in\F,\;0\le s_{ij}\le 1\;
\forall i\in\cS,\forall j\in\F_i.
\ea
\eeq
\end{theorem}
\proof{Proof.}
Without capacity constraints, it is easy to see that
the optimal objective value of the scenario problem
is given by
\beq\label{eqn:max-u*y-max-D*y}
\ba
&\sum_{i\in\cS}\big(\max_{j\in\F_i}u_{ij}y_j\big)  \big(\max_{k\in\F_i}D^{\omega}_{ik}y_k\big) 
=\sum_{i\in\cS} \max_{j,k\in\F_i} u_{ij}D^{\omega}_{ij} y_jy_k
=\sum_{i\in\cS}\max_{j\in\F_i}u_{ij}D^{\omega}_{ij}y_j,
\ea
\eeq
where we use the assumption that the utility gain and demand are consistent.
The problem \eqref{opt:FL-DRO} becomes the following:
\beq\label{opt:DRO-FL1}
\ba
&\underset{y}{\max}\;\underset{\mu}{\min}\; \sum_{\omega\in\Omega}\mu^\omega\sum_{i\in\cS} \max_{j\in\F_i} u_{ij}D^{\omega}_{ij} y_j \\
&\textrm{ s.t. }\; \sum_{j\in\F} b_jy_j\le B, \\
&\qquad\;\sum_{\omega\in\Omega}|\mu^{\omega}-\mu^{\omega}_0|\le d, \\
&\qquad\; y_j\in\{0,1\}\;\forall j\in\F.
\ea
\eeq
The terms involved in the inner minimization problem of \eqref{opt:DRO-FL1}
can be written as
\beq
\ba
&\underset{\mu}{\min}\;\sum_{\omega\in\Omega}\mu^\omega 
\Big(\sum_{i\in\cS} \max_{j\in\F_i} u_{ij}D^{\omega}_{ij} y_j\Big)  \\
&\textrm{ s.t. }\;\sum_{\omega\in\Omega} \rho^{\omega}\le d,  \\
&\qquad\; \mu^{\omega}-\mu^{\omega}_0\le\rho^{\omega} \qquad\forall \omega\in\Omega, \\
&\qquad\; \mu^{\omega}_0-\mu^{\omega}\le\rho^{\omega} \qquad\forall \omega\in\Omega, \\
&\qquad\; \sum_{\omega\in\Omega} \mu^{\omega}=1, \\
&\qquad\; \mu^{\omega}\ge 0,\;\rho^{\omega}\ge 0\; \forall \omega\in\Omega.
\ea
\eeq
Taking the dual of the above linear program with respect to the probability measure $\mu$, 
we obtain the following inner problem:
\bdm
\ba
&\max\; \sum_{\omega\in\Omega}\mu^{\omega}_0(\alpha^{\omega}-\beta^{\omega})-\lambda d-\gamma \\
&\textrm{ s.t. }\; \sum_{i\in\cS} \max_{j\in\F_i} u_{ij}D^{\omega}_{ij} y_j +\alpha^{\omega} 
-\beta^{\omega}+\gamma \ge 0 \qquad \forall\omega\in\Omega, \\
&\qquad \lambda-\alpha^{\omega}-\beta^{\omega}\ge 0 \qquad \forall\omega\in\Omega, \\
&\qquad \lambda\ge0,\;\alpha^{\omega}\ge0,\;\beta^{\omega}\ge0,\;\gamma\in\mathbb{R}.
\ea
\edm
To linearize the term $\max_{j\in\F_i} u_{ij}D^{\omega}_{ij} y_j$,
we can introduce binary indicator variables $s_{ij}$, and reformulate
the first constraint as follows
\bdm
\ba
&\sum_{i\in\cS}\sum_{j\in\F_i}u_{ij}D^{\omega}_{ij}s_{ij}+\alpha^{\omega} 
-\beta^{\omega}+\gamma \ge 0, \\
&s_{ij}\le y_j, \; \sum_{j\in\F_i}s_{ij}=1, \; 0\le s_{ij}\le 1,
\ea
\edm 
where we implicitly use the consistency condition of utility and demand
which implies $\textrm{argmax}_{j\in\F_i} u_{ij}D^{\omega}_{ij} y_j$ is scenario
independent and hence the variable $s_{ij}$.
Incorporate with the outer maximization yields the reformulation \eqref{opt:FL-DRO-uncap}.\Halmos
\endproof

\section{A Statistical Model for Utilities and Demand Estimation}
We establish a regression model based on the maximum attraction principle
for estimating utility value $u_{ij}$ and demand $D_{ij}$ that are input parameters to the service center
location model \eqref{opt:FL-DRO}. As a by-product, this regression model can also be
used to estimate the size parameter $d$ of the ambiguity set \eqref{eqn:amb-set}.
The regression model requires samples of survey among potential customers from a
site $i\in\cS$ on their ratings of multiple candidate service center locations.
Before formulating the regression model, we first introduce how the samples of survey
are collected for fitting the regression model. Suppose $N$ residents have 
been randomly selected from the site $i$, and they are viewed as potential customers
of service centers under planning. Each of them is asked to give a score in the range
$\{0,1,\ldots,q\}$ to each candidate service center location in $\F_i$. The score that
a customer is assigned to a location $j\in\F_i$ is taken as the potential utility 
gained by the customer if going to the service center at $j$. Score value 0
means the customer is unwilling to get service from the corresponding location.
Since the regression model structure is identical for each $i\in\cS$, 
we omit the customer site index $i$ and re-write $\F_i$ as $\cG$ in the following
modeling and analysis. Suppose the indices of location in $\cG$ are labeled
as $\cG:=\{1,2,\ldots,g\}$, the $N$ residents are labeled as $\{1,2,\ldots,N\}$,
and the score assigned to a candidate location $j\in\cG$ by the resident $k$
is denoted as $a_{kj}$.

The maximum attraction principle described in \eqref{eqn:max-attr}
may not be satisfied in reality. But it is possible to establish a utility-demand
estimation model that approximately meets the maximum attraction principle.
For example, we can verify whether the samples from survey satisfy the maximum
attraction principle by grouping the customers who have taken the survey
as follows:
\beq\label{def:Vj}
V_j=\Set*{k\in[N]}{a_{kj}\ge 1} \qquad\forall j\in\cG,
\eeq
where $[N]:=\{1,\ldots,N\}$. In words, $V_j$ is the set of residents who 
are willing to go to a service center located at $j$. According to \eqref{eqn:max-attr}, 
the samples satisfy the maximum attraction principle exactly
if there exists a permutation $\sigma$ on the indices in $\cG$
such that the following inclusive condition holds:
\beq
V_{\sigma(1)}\subseteq V_{\sigma(2)}\subseteq \ldots \subseteq V_{\sigma(g)}.
\eeq
A logic behind the maximum attraction principle is that
if two candidate locations have similar features, then if a customer
is attracted by one location, the customer should also be attracted
by the other location with high chance. In this case,
the difference in the attractability of the two candidate locations viewed by
the customer is reflected in the score assigned to the two locations
by the customer. The scores in this case are both non-zero indicating
that the customer is willing to visit any of them.
Only in the case that two candidate locations have some 
substantial differences (i.e., one is too far away from the customer site
or one is located at a bad community), a customer will be 
willing to visit one location (assigning a non-zero score to it)
while unwilling to visit the other one (assigning a zero score to it).

\subsection{An inclusive chain representation of the maximum attraction principle}
We establish an inclusive chain model to represent the set-level 
realization of the maximum attraction principle. First we build the
subsets \eqref{def:Vj} using collected samples, and the we sort
the indices in $\cG$ such that 
$|V_{\sigma(1)}|\le|V_{\sigma(2)}|\le\ldots\le|V_{\sigma(g)}|$,
where $\sigma$ is a permutation on $\cG$ that makes this
condition hold. For clarity we assume that 
$|V_1|\le|V_2|\le\ldots\le|V_g|$ without loss of generality.
The inclusive chain model assumes that the score assigned to
the candidate locations in $\cG$ by a customer from a fixed 
location should match with one of the following patterns:
\beq\label{eqn:valid-pattern}
[0,\ldots,0,a^{(i)}_i,\ldots,a^{(i)}_g] \qquad\textrm{for } i\in\{1,\ldots,g+1\},
\eeq
where in the $i$-th pattern, first $i-1$ scores are all zero and the
remaining $g-i+1$ scores are all non-zero. 
By convention, the 
$(g+1)$-th pattern is just $[0,\ldots,0]$. In a probabilistic flavor,
the inclusive chain model (ICM) for the score vector 
can be formally established as a
statistical model presented as
\beq\label{eqn:ICM}
\ba
&\xi\overset{\mathcal{D}}{=\joinrel=} \sum^{g+1}_{i=1}\bs{1}\{B=i\}[0,\ldots,0,Q^{(i)}_i,\ldots,Q^{(i)}_g], \\
&P(Q^{(i)}_j=r)=\pi^{(i)}_{jr} \qquad\forall i\in\{1,\ldots,g\},\;\forall j\in\{i,\ldots,g\},\;\forall r\in\{1,\ldots,q\}, \\
&P(B=i)=p_i\;\forall i\in\{1,\ldots,g+1\},\; \sum^{g+1}_{i=1}p_i=1, \\
&\sum^q_{r=1}\pi^{(i)}_{jr}=1\quad\forall i\in\{1,\ldots,g\},\;\forall j\in\{i,\ldots,g\},
\ea\tag{ICM}
\eeq 
where $\xi$ is the random score vector, $B$ is an indicator random
variable that selects the pattern matching with the score vector,
$p_j$ is the probability
of score vector matching with the $j$-th pattern,
and $\Set*{\pi^{(i)}_{jr}}{r=1,\ldots,q}$ gives the probability distribution
of the random score $Q^{(i)}_j$. Note that $p_i$ and $\pi^{(i)}_{jr}$
are model parameters that can be determined by fitting the inclusive
chain model with collected samples. The following proposition
connects the inclusive chain model with the set-level maximum attraction principle.
\begin{proposition}
If the score vector of every customer taking the survey follows the inclusive chain model,
then the set-level maximum attraction principle 
$V_1\subseteq V_2\subseteq \ldots \subseteq V_g$ 
is satisfied by the collection of samples almost surely.
\end{proposition}
\proof{Proof.}
We prove it contradiction. Suppose there exist subsets $V_i$
and $V_j$ (with $i<j$) such that $V_i\setminus V_j\neq\emptyset$
with some positive probability. Suppose $k$ is the customer 
who is willing to visit $V_i$ but not $V_j$, and let $\xi^k$ be the score
vector of this customer. The definition of $V_i$ implies that
$\xi_i\ge 1$ but $\xi^k\notin V_j$ implies $\xi^k_j=0$. 
On the other side, for every pattern
vector $\zeta$ we should have $\bs{1}\{\zeta_i>0\}\le\bs{1}\{\zeta_j>0\}$
almost surely, which implies that $\bs{1}\{\xi^k_i>0\}\le\bs{1}\{\xi^k_j>0\}$
almost surely. But this contradicts to $\xi^k_i\ge 1$ and $\xi^k_j=0$. \Halmos
\endproof

\subsection{An adjusted inclusive chain model for incorporating defective score vectors}
\label{sec:aicm}
We consider the problem of fitting the inclusive chain model with the collected samples (score vectors).
The first step is to determine the rank (level of attracability) 
of candidate locations $\cG$ involved in the inclusive chain.
An empirical method by simply sorting the cardinality of $V_i$
can be used to achieve this. The probability guarantee of this method which will be discussed in a moment.
For now, assume that the order has been identified, and suppose 
$|V_1|\le |V_2|\le \ldots |V_g|$ without loss of generality. So empirically,
this implies that the rank is $1,2,\ldots,g$ sorted by level of attractability 
from low to high. Based on this information, we can establish an inclusive
chain model as \eqref{eqn:ICM}, but this model is not capable 
to handle score vectors in which there is at least one zero-value entry between
two non-zero value entries, i.e., $[0,\ldots,0,5,7,0,0,7,6,9]$. 
We call a score vector \emph{defective} if it can not match with any pattern
vector in \eqref{eqn:valid-pattern} valid for the inclusive chain model.
Note that to convert the defective score vector $[0,\ldots,0,5,7,0,0,7,6,9]$ 
into a qualified pattern 
one can either add two 1 score entries to the middle 
($[0,\ldots,0,5,7,0,0,7,6,9]\to[0,\ldots,0,5,7,1,1,7,6,9]$)
or remove the first two nonzero entries ($[0,\ldots,0,5,7,0,0,7,6,9]$$\to$$[0,\ldots,0,7,6,9]$).
The first approach has added a total 2 units of score while the second
one has removed a total 12 units of score, and it can be verified that the first
approach is the one that makes the minimum total score change among all
possible changes that can convert this score vector into a qualified pattern.
Based on this rule, we can define the defective score units for an arbitrary score vector
as follows.
\begin{definition}\label{def:defect-score-qt}
For a given arbitrary score vector $v$, the \emph{defective score quantity}
is given as 
\beq
\min\;\norm{v-w}_1\; \textrm{ s.t. } w \textrm{ is a qualified pattern of the inclusive chain model}.
\eeq
\end{definition}
\begin{observation}
The defective score quantity has an upper bound $g-1$. 
\end{observation}
Let $H_s$ be the event that the defective score quantity of score vector is $s$,
and $m_s$ is the probability of this event. After introducing these additional
parameters, the inclusive chain model can then be adjusted to incorporate
all defective score vectors as follows
\beq\label{eqn:AICM}
\ba
&\E[\xi|H_0]\overset{\mathcal{D}}{=\joinrel=} \sum^{g+1}_{i=1}\bs{1}\{B=i\}[0,\ldots,0,Q^{(i)}_i,\ldots,Q^{(i)}_g], \\
&P(Q^{(i)}_j=r)=\pi^{(i)}_{jr} \qquad\forall i\in\{1,\ldots,g\},\;\forall j\in\{i,\ldots,g\},\;\forall r\in\{1,\ldots,q\}, \\
&P(B=i)=p_i\;\forall i\in\{1,\ldots,g+1\},\; \sum^{g+1}_{i=1}p_i=1, \\
&\sum^q_{r=1}\pi^{(i)}_{jr}=1\quad\forall i\in\{1,\ldots,g\},\;\forall j\in\{i,\ldots,g\}, \\
&P(\xi\in H_s)=m_s \quad\forall s\in\{0,1,\ldots,g-1\},\; \sum^{g-1}_{s=0}m_s=1.
\ea\tag{AICM}
\eeq 
The above model is referred as the adjusted inclusive chain model \eqref{eqn:AICM}. 

\subsubsection{A theoretical probability guarantee of identifying the order of the inclusive chain model}
As mentioned at the beginning of Section~\ref{sec:aicm}, 
a simple sorting method can be used to identify the rank 
of the candidate locations in $\cG$ for the AICM.
We sort the cardinality of subsets $V^{\prime}_j$ ($j\in\cG$) defined as 
\bdm
V^{\prime}_j=\Set*{k\in[N]}{a_{kj}\ge 1},
\edm
where $a_k$ is the score vector of the $k$-th customer participated in the survey,
and $a_{kj}$ is its $j$-th element. The following proposition provides a probability
guarantee on identifying the order of inclusive chain model given that samples obey
the AICM.
\begin{proposition}
Suppose the set of $N$ samples collected from survey follow the AICM,
Let the set of candidate locations be $\cG=\{1,\ldots,g\}$.
Let $E$ be the event that
the rank of inclusive chain model can be successfully identified by sorting the
cardinality of $V^{\prime}_j$ ($j\in\cG$). Let $p^*=\min_{j\in\cG}p_j$.
If $p^*>(1-m_0)/m_0$, then 
\beq
P(E)\ge 1-(|\cG|-1)\emph{exp}\Big(-\frac{1}{2}N[m_0(1+p^*)-1]^2\Big).
\eeq
\end{proposition}
\proof{Proof.}
We first show that with the given assumptions,
the rank of the inclusive chain model can be identified almost surely as the total number
of samples $N$ goes to infinity. In the analysis, we let $\sigma$ be the permutation that
gives the rank of locations, i.e., $\sigma(1)<\ldots<\sigma(g)$ is the rank of locations from low to high. 
Let $N_j$ be the number of samples that match with
the $j$-th pattern $[0,\ldots,0,Q^{(j)}_{\sigma(j)},\ldots,Q^{(j)}_{\sigma(g)}]$.
Based on the probability setting of AICM, and the law of large numbers we should have the following asymptotic
relation hold:
\beq
\frac{N_j}{N}\to m_0p_j \qquad a.s. \textrm{ as } N\to\infty \qquad\forall j\in\cG. 
\eeq
This implies that $N_j>0$ almost surely as $N\to\infty$. 
By the definition of $V^{\prime}_j$, we should have the following bounds
for every $|V^{\prime}_j|$
\beq
\sum^{\sigma(j)}_{k=1}N_k\le |V^{\prime}_j|\le \sum^{\sigma(j)}_{k=1}N_k+|\cup^{g-1}_{s=1}H_s|  \qquad\forall j\in\cG.
\eeq
For any $i,j$ satisfying $\sigma(i)<\sigma(j)$, it suffices to show that 
$|V^{\prime}_{\sigma(i)}|<|V^{\prime}_{\sigma(j)}|$ as $N\to\infty$. Indeed we have
\bdm
\ba
&\frac{|V^{\prime}_{\sigma(j)}|-|V^{\prime}_{\sigma(i)}|}{N}\ge
\frac{1}{N}\Big(\sum^{\sigma(j)}_{k=1}N_k - \sum^{\sigma(i)}_{k=1}N_k - |\cup^{g-1}_{s=1}H_s| \Big) \\
&=\frac{1}{N}\Big(\sum^{\sigma(j)}_{k=\sigma(i)+1}N_k - |\cup^{g-1}_{s=1}H_s| \Big)
\to \sum^{\sigma(j)}_{k=\sigma(i)+1}mp_{k}-(1-m_0) \quad \textrm{ as }N\to\infty,
\ea
\edm
and $\sum^{\sigma(j)}_{k=\sigma(i)+1}mp_{k}-(1-m_0)\ge mp^*-(1-m_0)>0$ by assumption.
Therefore the rank of the inclusive chain can be 
identified almost surely as $N$ goes to infinity. 
Now we prove the non-asymptotic result. 
Note the probability $P(E)$ can be lower bounded as follows
\bdm
\ba
&P(E)=P\Big(\bigcap^{g-1}_{i=1}\big\{|V^{\prime}_{\sigma(i+1)}|-|V^{\prime}_{\sigma(i)}|>0\big\} \Big) 
=1-P\Big(\bigcup^{g-1}_{i=1}\big\{|V^{\prime}_{\sigma(i+1)}|-|V^{\prime}_{\sigma(i)}|\le 0\big\} \Big) \\
&\ge 1-\sum^{g-1}_{i=1}P\Big(|V^{\prime}_{\sigma(i+1)}|-|V^{\prime}_{\sigma(i)}|\le 0\Big) 
\ge 1-\sum^{g-1}_{i=1}P\Big( N_{\sigma(i+1)} - |\cup^{g-1}_{s=1}H_s|\le 0 \Big).
\ea
\edm
Since we have $N_{\sigma(i+1)} - |\cup^{g-1}_{s=1}H_s|=\sum^N_{k=1} I_k$,
where $\{I_k\}^N_{k=1}$ are i.i.d. random variables satisfying $P(I_k=1)=m_0p_{i+1}$,
$P(I_k=0)=m_0(1-p_{i+1})$ and $P(I_k=-1)=1-m_0$, which gives 
$\E[I_k]=m_0(1+p_{i+1})-1>0$ by assumption. 
Let $\bar{I}=\frac{1}{N}\sum^N_{k=1}I_k$
Applying Hoeffding's inequality gives the following bound
\bdm
\ba
&P\Big( N_{\sigma(i+1)} - |\cup^{g-1}_{s=1}H_s|\le 0 \Big)=P\Big(\sum^N_{k=1} I_k\le 0 \Big) \\
&=P\Big(\bar{I}-\E[\bar{I}]\le 1-m_0(1+p_{i+1}) \Big)\le\exp\Big(-\frac{1}{2}N[m_0(1+p_{i+1})-1]^2\Big)\\
&\le \exp\Big(-\frac{1}{2}N[m_0(1+p^*)-1]^2\Big).
\ea
\edm
Then it follows that $P(E)\ge 1-(|\cG|-1)\exp\Big(-\frac{1}{2}N[m_0(1+p^*)-1]^2\Big)$. \Halmos
\endproof

\subsubsection{Model fitting and parameter estimation}
To avoid redundant notations in the following analysis, 
we assume that candidate locations are ranked as $1<\ldots<g$
based on attractability. 
To fit the AICM with samples, a likelihood function
and some regularity terms need to be established. Note that without consideration
of any regularity, all model parameters can be estimated naturally using a frequency
count approach based on the given samples. For example, an empirical estimation of  
$p_j$ can be $\hat{p}_j=N_j/N^{\prime}$, and an empirical estimation of $\pi^{(i)}_{jr}$
can be $\hat{\pi}^{(i)}_{jr}=\frac{N^{(i)}_{jr}}{\sum^q_{r=1}N^{(i)}_{jr}}$, where
$N^{(i)}_{jr}$ is the number of samples that have $Q^{(i)}_{j}=r$. However, this 
estimation can lead to overfitting of new samples. We now focus on establishing
a likelihood function and introduce some regularity terms to prevent overfitting.
First, the parameter $m_s$ can be estimated as 
$\wh{m}_s=\frac{\#\{k\in[N]:\;\xi_k\in H_s\}}{N}$ for all $s\in\{1,\ldots,g-1\}$.
All defective samples can be further used to estimate other parameters 
with more information. If a sample $\xi_k$ is defective, it can be converted into 
a qualified score vector $\xi^{\prime}_k$ by imposing a 
defective score quantity (Definition~\ref{def:defect-score-qt}) 
to $\xi_k$, and hence the original samples 
$\{\xi_k\}^N_{k=1}$ can then be converted into a set of qualified score
vectors $\{\xi^{\prime}_k\}^N_{k=1}$. We can evaluate the probability
of generating $\xi^{\prime}_k$ by the model using the model parameters,
and denote this probability as $P(\xi^{\prime}_k|p,\pi)$. Then the (logarithmic) likelihood
function of generating $\{\xi^{\prime}_k\}^N_{k=1}$ can then be formulated as
\beq
\mathcal{L}(p,\pi)=\frac{1}{N}\sum^N_{k=1}\log P(\xi^{\prime}_k|p,\pi).
\eeq
Note that maximizing the above likelihood function with only probability-normalization constraints
leads to a frequency estimation of all parameters. 
Some regularity terms can be added to balance this trend. 
First, for score vectors in the same pattern, parameters that lead to smaller variance of score
on a given candidate location are preferred, and hence the following regularity terms
can be added (with some regularity coefficients):
\bdm
\ba
R^1_{ij}(\pi):=\var[Q^{(i)}_j]=\sum^q_{r=1}\Big(r^2\pi^{(i)}_{jr} \Big)-\Big(\sum^q_{r=1}r\pi^{(i)}_{jr}\Big)^2 \\
\qquad\forall i\in\{1,\ldots,g\},\;\forall j\in\{i,\ldots,g\}.
\ea
\edm
Second, consider any two consecutive locations $j$ and $j+1$ in the inclusive chain.
The two locations are both included in patterns $H_1,\ldots,H_j$. In the pattern $H_i$ 
($1\le i\le j$), the mean score difference between the two locations is given by
\bdm
w^i_j=\E[Q^{(i)}_{j+1}-Q^{(i)}_j]= \sum^q_{r=1} r(\pi^{(i)}_{j+1,r}-\pi^{(i)}_{j,r}).
\edm
For every pattern that contains locations $j$ and $j+1$, we expect the above
mean score difference to be similar. To enforce this similarity, we can add the following
regularity terms (with some regularity coefficients):
\bdm
R^2_j(\pi):=\frac{1}{j}\sum^{j}_{i=1}(w^i_j)^2-\frac{1}{j^2}\Big(\sum^j_{i=1}w^i_j\Big)^2 \qquad\forall j\in\{1,\ldots,g-1\}.
\edm
After incorporating the above regularity terms, the parameters can be estimated
by solving the following constrained minimization problem of the loss function.
\beq\label{opt:model-fit}
\begin{aligned}
&\min\; L(p,\pi):=\frac{2\lambda_1}{g(g+1)}\sum^g_{i=1}\sum^g_{j=i}R^1_{ij}(\pi)
+\frac{\lambda_2}{g-1}\sum^{g-1}_{j=1}R^2_j(\pi)
-\frac{1}{N}\sum^N_{k=1}\log P(\xi^{\prime}_k|p,\pi) \\
&\textrm{ s.t. } \sum^q_{r=1}\pi^{(i)}_{jr}=1 \qquad \forall i\in\{1,\ldots,g\},\;\forall j\in\{i,\ldots,g\},  \\
&\qquad \sum^{g+1}_{j=1}p_j=1,
\end{aligned}
\eeq
where $\lambda_1$ and $\lambda_2$ are regularity coefficients that can be tuned 
using a cross-validation approach. Note that we add a negative sign to the 
logarithmic likelihood function when it is incorporated into the loss function $L(p,\pi)$.
Notice that \eqref{opt:model-fit} is a nonconvex optimization
problem with continuous parameters and simple constraints, 
and it can be solved to local optimality very efficiently using
numerical optimization solvers. 

\subsection{Generation of utility and demand samples guided by AICM}
\label{sec:gen-dmd-amb-set}
Once the regression problem \eqref{opt:model-fit} is solved to give an estimation of the AICM parameters, 
the AICM can be used to estimate the (mean) utility value, guide the generation 
of demand samples.
First, the utility value $u_{\cdot,j}$ (the customer site index is ignored since the 
AICM is established for any given customer site) can be estimated as follows:
\beq
\hat{u}_{\cdot,j}=\sum^j_{i=1}\hat{p}_i\Big(\sum^q_{r=1}r\hat{\pi}^{(i)}_{jr}\Big),
\eeq
where $\hat{p}$ and $\hat{\pi}$ are the estimation of AICM parameters.
Generating a number of joint demand samples and use them to construct
the ambiguity set can be a challenging problem in practice. There are two major obstacles
in achieving this. 
(a) Scaling factor: since only a limited amount of survey samples are collected, how
to scale the demand estimation based on the samples to roughly estimate the 
total demand from a customer site? 
(b) Demand catenation: how to catenate the estimated demand 
samples for each customer site to get joint demand vectors (for all customer sites) 
as samples that can be directly input to the ambiguity set?
We try to give a few guidelines for dealing with these problems in this paper.

For (a), the total potential demand 
from a customer site could be estimated via a different channel
which is not the focus of this paper. Some investigation on the
demand estimation and especially the healthcare service demand
estimation is conducted in \citep{2005-dmd-est-heterog-consumers,
2001-imp-healthcare-dmd-estim,
2017-est-dmd-suppl-healthcare-serv,2008-opt-comm-health-center-loc}.
The survey of scoring on candidate locations are only
used to quantify the attractability of each candidate location
that can lead to different estimated demand for different candidate
locations when the total potential demand is given. 
In the next paragraph for addressing (b), we provide a estimation
of demand associated with each candidate location based on
assuming the total potential demand is $N$ (the number of collected samples).
In the case that we have a separate estimation $A$ (can be a random variable) 
of total potential demand the scaling factor can be defined as $\frac{A}{N}$,
and it can be used to scale up estimated demand associated with different 
candidate locations.
In the following discussion for (b), we assume that the demand samples
generated (using Algorithm~\ref{alg:gen-sp-Di}-\ref{alg:gen-defect-sp}) 
will be multiplied by the scaling factor (could be customer site dependent),
and omit this factor in these algorithms. 

For (b), one approach is to assume
the distributions of demand $D_i:=\Set*{D_{ij}}{j\in\F_i}$ for different customer sites $i$ are 
independent. In this case we can generate demand
samples for each $i$ and combine them (in an outer-product manner)
to construct samples for the random vector $D:=\Set*{D_{ij}}{i\in\cS,\;j\in\F_i}$.
For example, if we have constructed $|\Omega_i|$ samples for  
$D_i$ which are denoted as $\Set*{D^{\omega}_i}{\omega\in\Omega_i}$,
then we can get samples for $D$ as 
$\Set*{[D^{\omega_i}_i:\;i\in\cS]}{\omega_i\in\Omega_i\;\forall i\in\cS}$.
In this manner, the number of samples for $D$ can be as large as $\prod_{i\in\cS}|\Omega_i|$.
To reduce the sample size for $D$, one can select a desired number $N$ of 
combinations of samples for each $D_i$ at random. Another way is to assume
the demand from different costumer sites are strongly correlated.
As a consequence, we can divide 
the samples for each $D_i$ into several bands according to the magnitude
of the demand. Then we select samples from the same band for different
$D_i$'s to construct a sample for $D$. A possible way of defining the bands
can be based on different intervals of quantiles. For example, one can define 
bands corresponding to 5 intervals of percentage $[0,0.2],\;[0.2,0.4],\ldots,[0.8,1.0]$
for each $D_i$. This way of constructing samples (scenarios) has been used
in \citep{2020-ventilator-alloc-covid19}. Now let us focus on how to generate samples for each $D_i$.
The pseudo code for generating a sample for $D_i$ is given in Algorithm~\ref{alg:gen-sp-Di}.
Note that by repeatedly call this algorithm, one can generate any given number of samples
for $D_i$. A pseudo-code for generating a sample for $D$ is given as Algorithm~\ref{alg:gen-defect-sp}.

\begin{algorithm}[H]
	{\footnotesize
	\caption{\footnotesize An algorithm for generating a mixed-defective sample for $D_i$ ($i\in\cS$)
	following AICM.}\label{alg:gen-sp-Di}
	\begin{algorithmic}[1]
	\State{Input data: score vector samples $Y_i=\{\widehat{\xi}^k\}^N_{k=1}$ that have been collected 
	from the survey at the customer site $i$.}
	\State{Determine the rank of ICM based on the score vectors in $Y_i$. Sort the 
	indices in $\F_i$ (from low to high) based on this rank.}
	\State{Based on this rank, partition $Y_i$ into qualified subset $Y^q_i$ of samples
		 and defective subset $Y^d_i$ of samples. Let $N^q=|Y^q_i|$, $N^d=|Y^d_i|$ and $\rho=N^q/(N^q+N^d)$.}
	\State{Solve \eqref{opt:model-fit} to get an optimal estimation of probability parameters
		denoted as $p^*,\pi^*$ and substitute $p^*,\pi^*$ into \eqref{eqn:ICM}.}
	\State{Generate $N$ i.i.d. score vectors denoted as $\{\xi^k\}^{N}_{k=1}$ such
	that each score vector is generated using the way given as follows: 
	with probability $\rho$ sample $\xi^k$ following the fitted \eqref{eqn:ICM},
	and with probability $1-\rho$ sample $\xi^k$ uniformly from $Y^d_i$.}
	\State{Let $\widehat{D}_{ij}=\textrm{card}\{k\;|\;\xi^k_j\ge 1\}$ for all $j\in\F_i$.} \label{lin:score-to-demand}
	\State{Return the sample $\widehat{D}_i=[\widehat{D}_{ij}:\;j\in\F_i]$.}
	\end{algorithmic}
	}
\end{algorithm}

\begin{algorithm}[H]
	{
	\caption{\footnotesize An algorithm for generating a mixed-defective sample for $D$ following AICM.}\label{alg:gen-defect-sp}
	\begin{algorithmic}[1]
	\State{Input data: (1) same input as Algorithm~\ref{alg:gen-sp-Di} for all $i\in\cS$. 
	   (2) sample size $N_i$ for every $i\in\cS$ and $N$.}
	 \State{Generate $N_i$ i.i.d. mixed-defective samples for every $i\in\cS$ using Algorithm~\ref{alg:gen-sp-Di}.}
	 \State{Catenate the mixed-defective demand samples $\{\widehat{D}^k_i\}^{N_i}_{k=1}$
	 over $i\in\cS$ (either in an independent approach or a strong correlated approach) 
	 and draw $N$ mixed-defective samples $\{\widehat{D}^k\}^{N}_{k=1}$ for $D$.}
	 \State{Return the mixed-defective samples $\{\widehat{D}^k\}^{N}_{k=1}$.}
	\end{algorithmic}
	}
\end{algorithm}

\section{Numerical Investigation}
The numerical investigation consists of two major parts:
the numerical analysis of the AICM and the computational performance
of the \eqref{opt:FL-DRO} model. In Section~\ref{sec:num-aicm}, 
we provide numerical analysis of the AICM, where we simulate
certain amount of score vector samples and use them to fit the AICM.
In Section~\ref{sec:num-FL-DRO}, we show the computational performance
of solving generated instances of the \eqref{opt:FL-DRO} model,
and provide sensitivity analysis of the robust optimal locations
with respect to the sample size and the choice of radius in the ambiguity set.

\subsection{Numerical analysis of the AICM}
\label{sec:num-aicm}
For the numerical analysis of AICM,
we focus on a specific customer site with 4 candidate service center
locations labeled as $L_1\sim L_4$ (the rank of attractability is unknown).
We consider 6 possible scores $\{0,1,\ldots,5\}$.   
We simulate $N$ i.i.d. score vectors with some defected ones
and use them to fit the AICM. The generation of each score vector
is based on the following procedures. First, we
simulate a 4-dimensional real vector $v$
following the Gaussian distribution $\mtc{N}(w,\Sigma)$
with $w=[1,3,4.5,2]$ corresponding to the mean score of 
$L_1\sim L_4$ respectively, and $\Sigma=\diag[0.8^2,1.5^2,1,1.5^2]$
as the covariance matrix. Next, each entry of $v$ is rounded to the closest integer
from $\{0,1,\ldots,5\}$. Finally, with probability 0.8 we accept this score vector,
while with probability 0.2 we randomly select an entry of $v$, 
set it to be zero and accept the resulting score vector. We follow the steps 
instructed in Section~\ref{sec:aicm} to fit the AICM. An essential step is to solve
the nonlinear optimization problem \eqref{opt:model-fit}. This problem is implemented
in Julia 1.5.2 \citep{julia} with optimization package JuMP dev \citep{JuMP} as the modeling interface
and Ipopt 3.13.2 \citep{ipopt} as the solver.
For simplicity, we let $\lambda_1=\lambda_2=\lambda$ in \eqref{opt:model-fit}.
We use a 4-fold cross validation approach to select an appropriate regularization
parameter $\lambda$. The parameter selection consists of two rounds. 
The first round is a rough selection, where the candidate $\lambda^{(1)}$ is selected 
from $\{0, 0.1, 0.2, $\ldots$, 5\}$ based on the out-of-sample performance.
Suppose $\lambda^{(1)}_*$ is the selected value after the first round.
In the second round, we select an optimal $\lambda^*$ from the set 
$\Set*{\lambda^{(1)}_*\pm 0.01k}{k=0,1,\ldots,10}$ based on the out-of-sample performance. 
The out-of-sample performance is measured by the averaged logarithm of likelihood.
The first-round selection results are given in Figure~\ref{fig:lam-sel} for 
three different sample sizes $N=50,100,200$. The second-round selection
results are given in Table~\ref{tab:lam-sel-sec-round}, from which we see that 
the optimal $\lambda$ is 2.98, 2.98 and 0.37 for $N=50,100,200$,
respectively. The curves in Figure~\ref{fig:lam-sel} indicate that for a small
sample size ($N=50$) a larger regularization
parameter lead to a better out-of-sample performance,
while for a larger sample size ($N=200$) the best out-of-sample performance
is achieved at a smaller regularization parameter. This behavior is intuitive
since training the model with a small sample size can lead to over fitting 
if the model parameters are not sufficiently regularized. When more samples
are available to train the model, the regularization can be less restrictive.

\begin{figure}
\centering
\begin{adjustbox}{minipage=\linewidth,scale=1}
\begin{subfigure}{.5\textwidth}
  \centering
  \includegraphics[width=.9\linewidth]{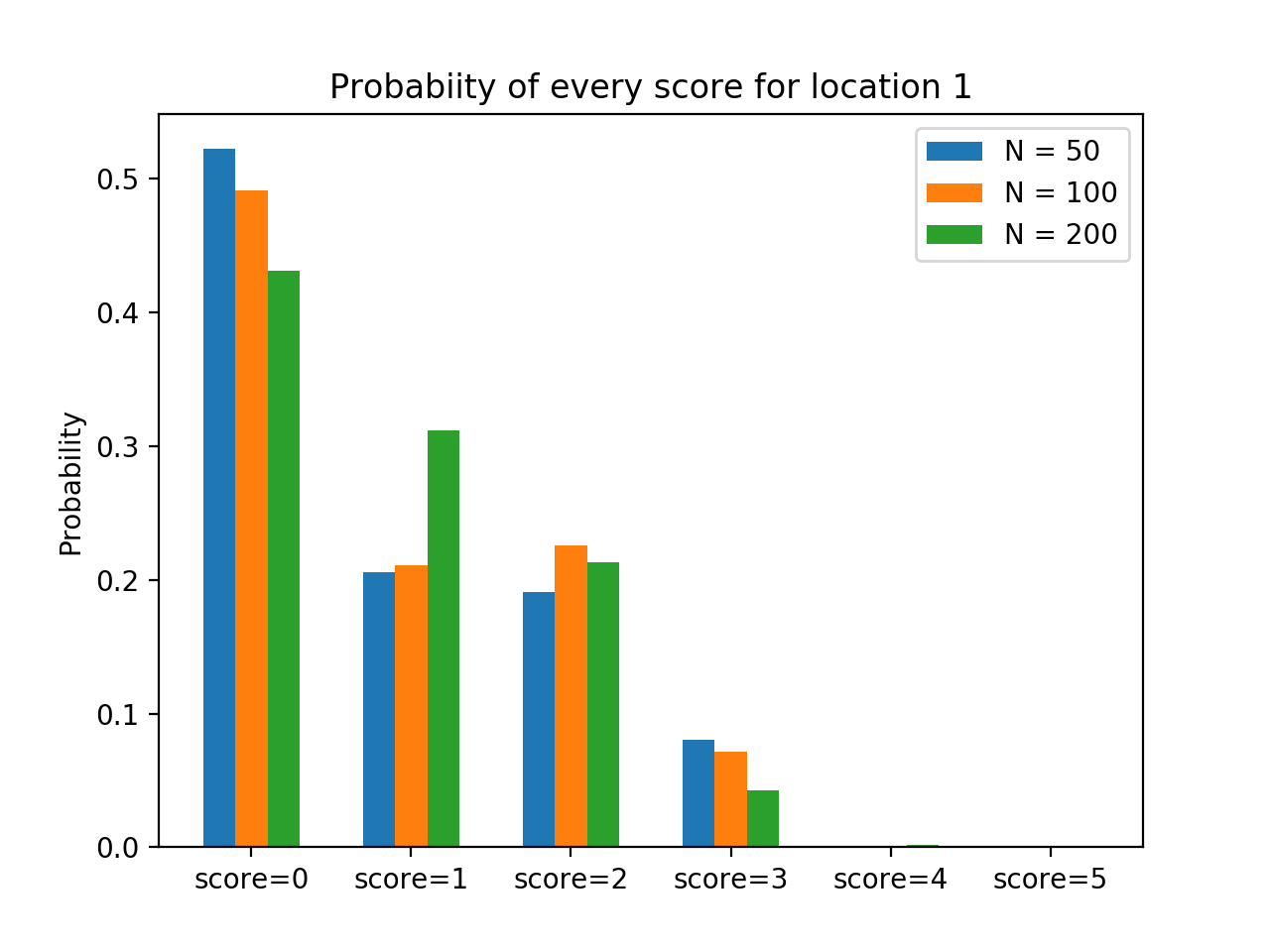}  
  \caption{Probability distribution on scores of location 1.}
\end{subfigure}
\begin{subfigure}{.5\textwidth}
  \centering
  \includegraphics[width=.9\linewidth]{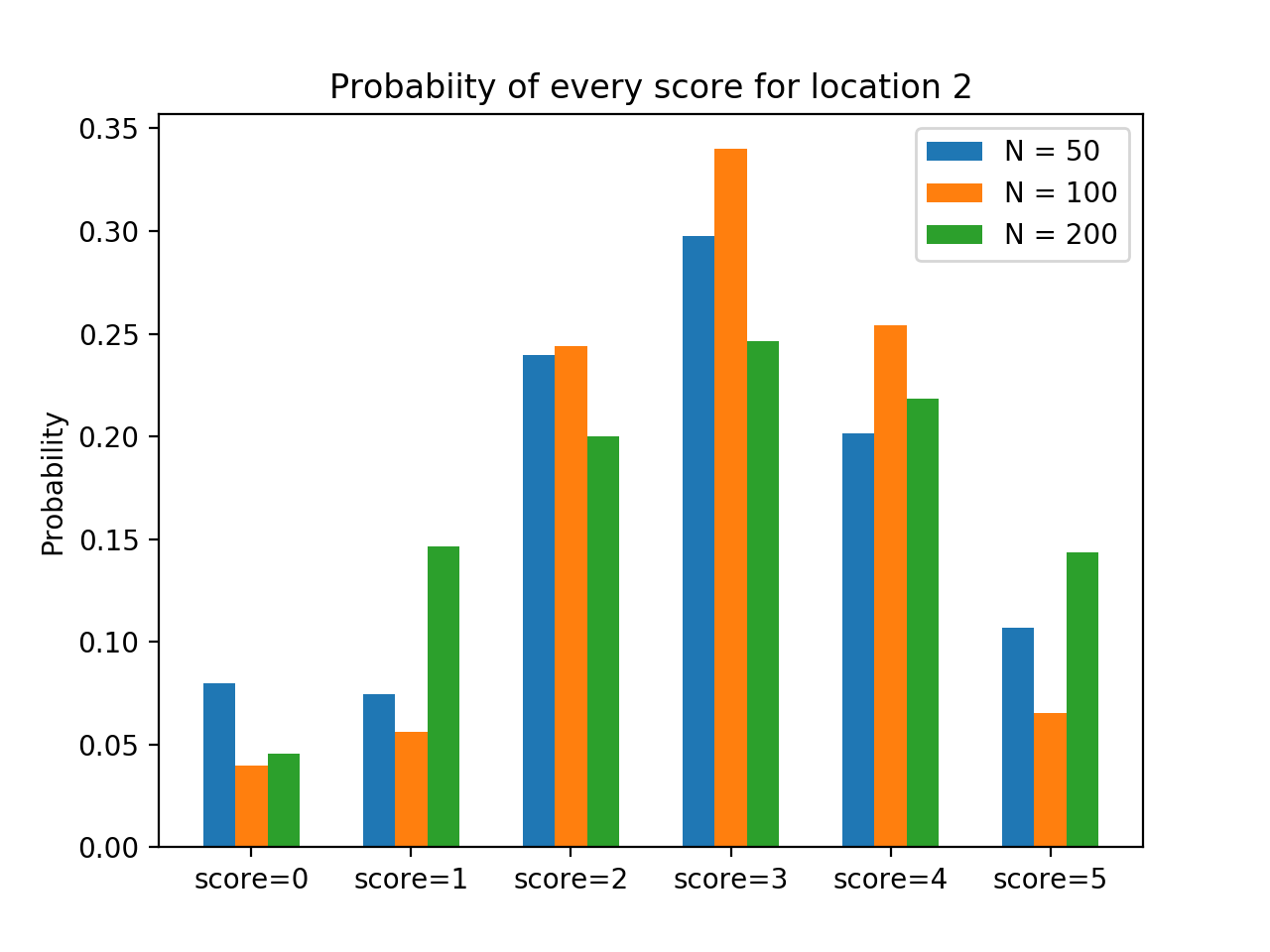}  
  \caption{Probability distribution on scores of location 2.}
\end{subfigure}
\begin{subfigure}{.5\textwidth}
  \centering
  \includegraphics[width=.9\linewidth]{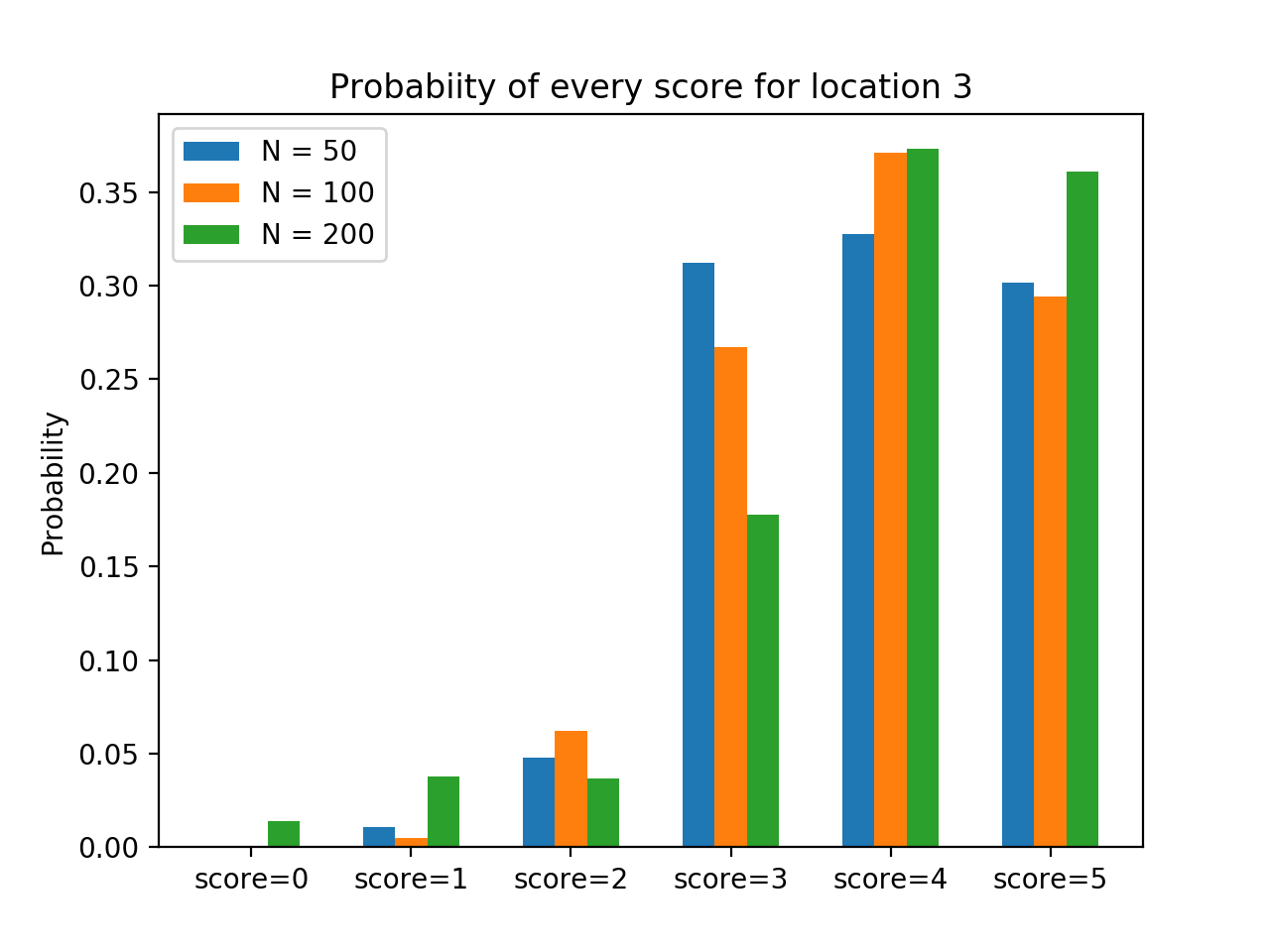}  
  \caption{Probability distribution on scores of location 3.}
\end{subfigure}
\begin{subfigure}{.5\textwidth}
  \centering
  \includegraphics[width=.9\linewidth]{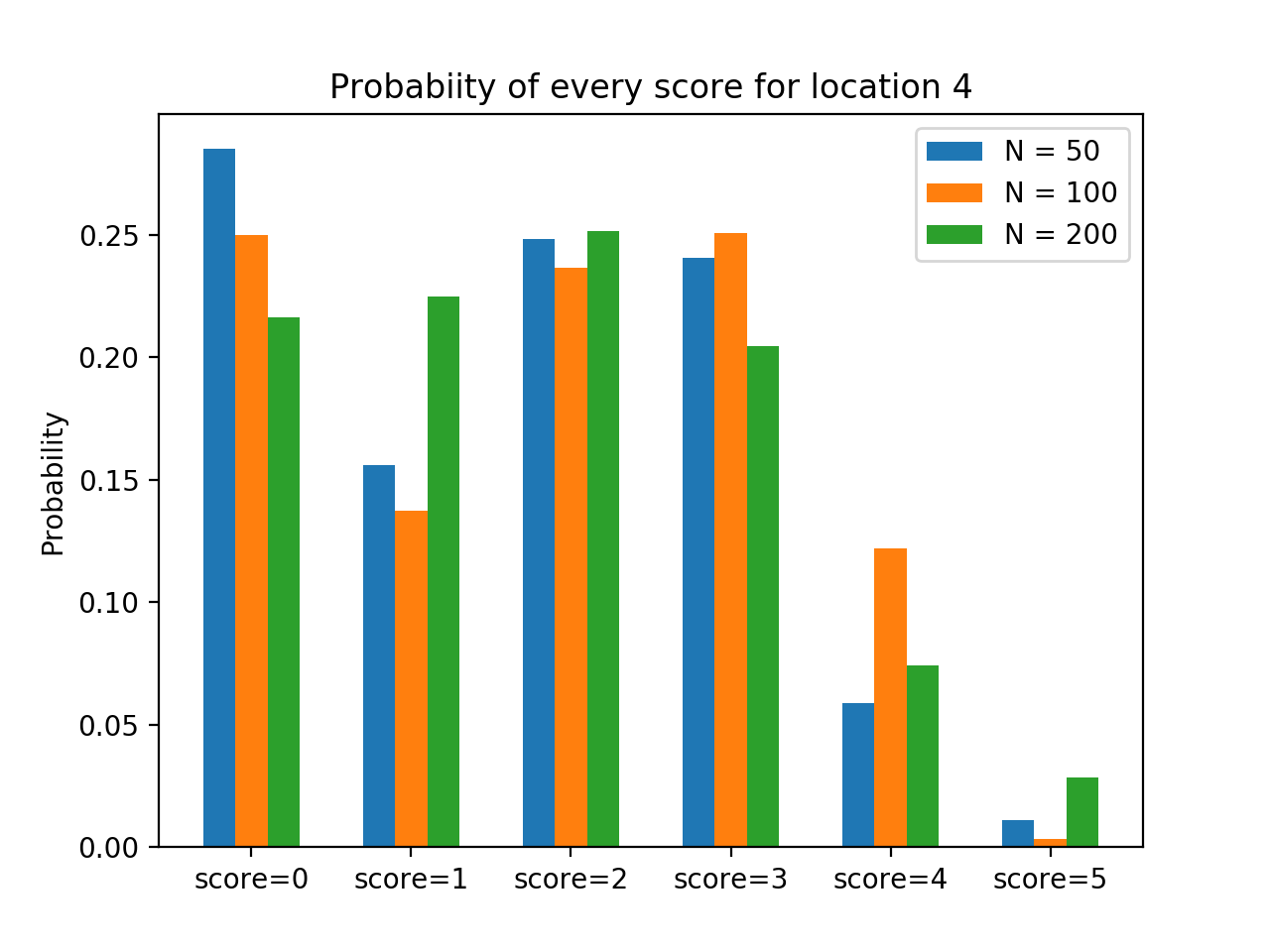}  
  \caption{Probability distribution on scores of location 4.}
\end{subfigure}
\end{adjustbox}
\caption{Probability distribution on scores for each location at different sample sizes.}
\label{fig:prob-distr-scores}
\end{figure}

Finally, we fit the AICM using $N$ samples and the optimal $\lambda$ parameter. 
The probability parameters $p$ and $\pi$ from AICM can be used to estimate the probability
distribution on scores for each candidate location. The results are given 
in Figure~\ref{fig:prob-distr-scores} for the 4 candidate locations. 
The probability distribution on score values reveal that the rank of attractability
from high to low is: Location 3 $>$  Location 2 $>$  Location 4 $>$  Location 1,
which matches with the underline stochastic model that generates the scores.
The mean scores corresponding to the four locations are given in Table~\ref{tab:mean-score}.

\begin{table}
\centering
\caption{Mean scores for each location at different sample sizes.}\label{tab:mean-score}
\begin{tabular}{r|rrrr}
\hline\hline
Samples	&	 Loc. 1	&	Loc. 2	&	Loc. 3	&	Loc. 4	\\
\hline
50	&	0.83	&	2.79	&	3.86	&	1.67	\\
100	&	0.88	&	2.91	&	3.89	&	1.87	\\
200	&	0.87	&	2.88	&	3.94	&	1.78	\\
\hline
\end{tabular}
\end{table}

\begin{table}
	\begin{minipage}{0.5\linewidth}
		\centering
		\includegraphics[width=.9\linewidth]{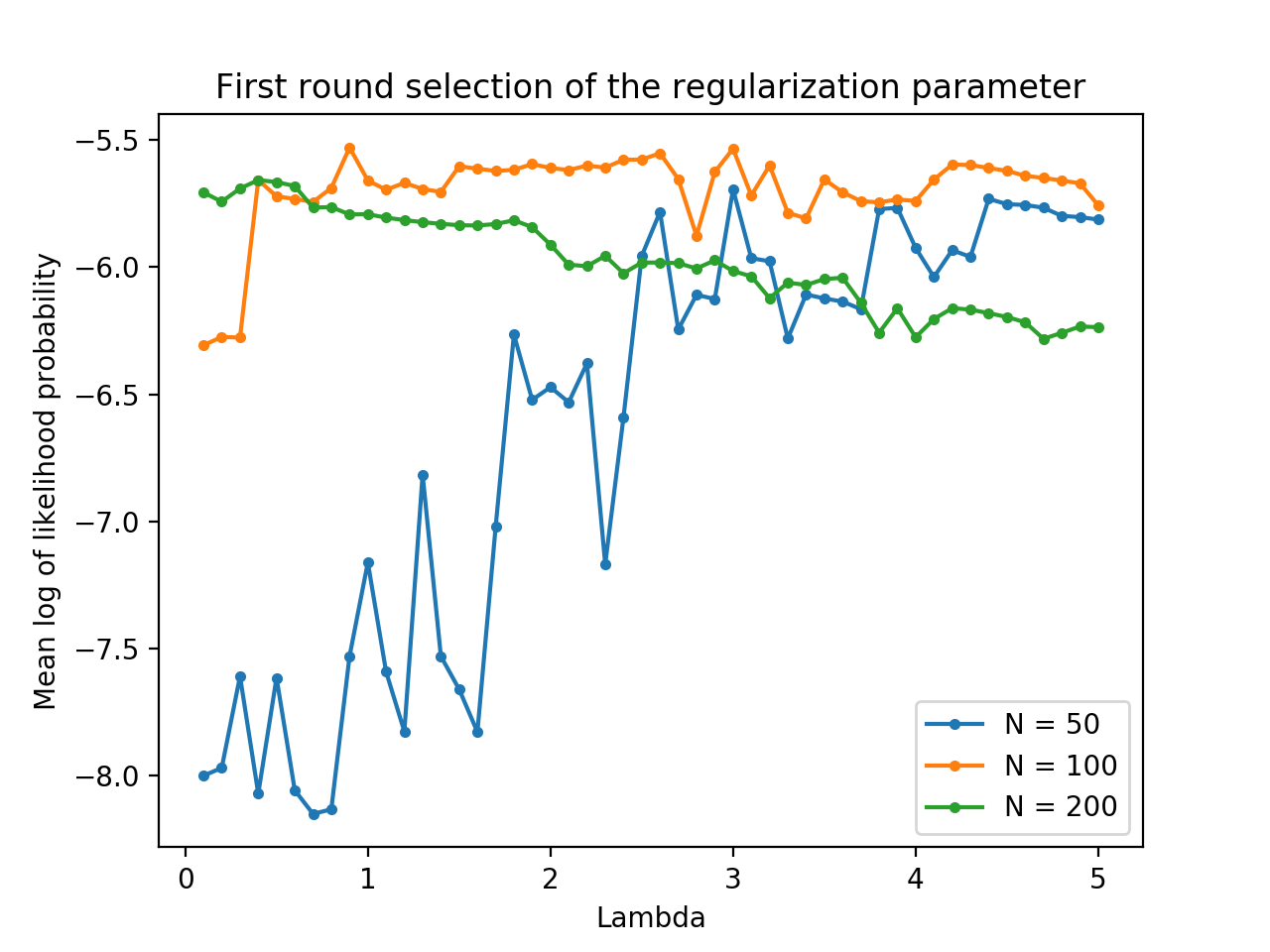}  
 		 \captionof{figure}{First round of $\lambda$ selection.}\label{fig:lam-sel}
	\end{minipage} 
	\begin{minipage}{0.3\linewidth}
		{
		\caption{Second round of $\lambda$ selection.}\label{tab:lam-sel-sec-round}
		\begin{tabular}{c|rrr}
		\hline\hline
		$N$	&  50  & 100  & 200	\\
		\hline
		range & [2.9,\;3.0]  &  [2.9,\;3.0]  & [0.3,\;0.4] \\
		$\lambda^*$ & 2.98 & 2.98 & 0.37 \\
		\hline
		\end{tabular}
		}
	\end{minipage}
\end{table}

\subsection{Computational performance of the \eqref{opt:FL-DRO} model}
\label{sec:num-FL-DRO}
We investigate the computational performance of solving numerical instances of
the \eqref{opt:FL-DRO} model. To generate numerical instances, we create a 
$100\times 100$ two-dimensional square as the map. The customer sties are points
inside the square, and the coordinates of each customer site are uniformly
drawn from the range $[0,100]^2$. We let $\F=\cS$, which means each customer
site can be a candidate location of a service center. 
We sort the set of distances $W=\Set*{d(i,j)}{\forall i\in\cS,\;\forall j\in\cS,\; i<j}$ 
from small to large, and let $d_c$ be the 5\% quantile of the sorted list.
For every customer site $i$, we select all customer sites (including $i$)
that are within $d_c$ distance from $i$ to form the set $\F_i$. 
We let the utility completely depend on the distance. Specifically, 
we set $u_{ii}=5$ for all $i\in\cS$. We also let the first
element in $W$ correspond to utility 5, the last element in $W$
correspond to utility 0.5, and any other elements correspond to utility values obtained
using a linear interpolation. For a sample $\omega$, the demand $D^{\omega}_{ij}$
is drawn from the normal distribution $\mtc{N}(120u_{ij}, (12u_{ij})^2)$.
The capacity of each candidate service center is set to be 1000.
We assume the cost of opening each service center is identical (an unit cost),
and the budget $B$ in this case is equivalent to the maximum number of service centers
that can be opened. We have generated 18 instances of \eqref{opt:FL-DRO}, and each
instance is characterized by the number $|\cS|$ of customer sites (also candidate locations)
and the budget $B$ (number of maximum service centers that can be allocated).
For each instance, we have tried three different sample-size options ($N=100,200,500$),
and the total variation distance is set to be 0.2 in these calculation.
The \eqref{opt:FL-DRO} model is implemented and solved using a 
Gurobi Python interface (Python 3.7 and Gurobi 9.0).
For every instance, the master problem and all scenario sub-problems are solved using a single 2.50GHz CPU. 
The computational results for the 12 instances are shown in Table~\ref{tab:comp-perform}.
It is shown that all instances have been solved to optimality
within 15 master iterations (see Algorithm~\ref{alg:cut-plane}).
The first 9 instances are solved within 120 seconds, and large
instances are solved within 1.5 hours where the majority of computational time
is spent on the scenario sub-problems.  
The range of the robust optimal objective is less than 0.7\% across
three different sample-size options for each instance, and we have also 
observed that the robust optimal solution remains the same for different
sample-size options.

\begin{figure}
\centering
\begin{subfigure}{.45\textwidth}
 \includegraphics[width=1\linewidth]{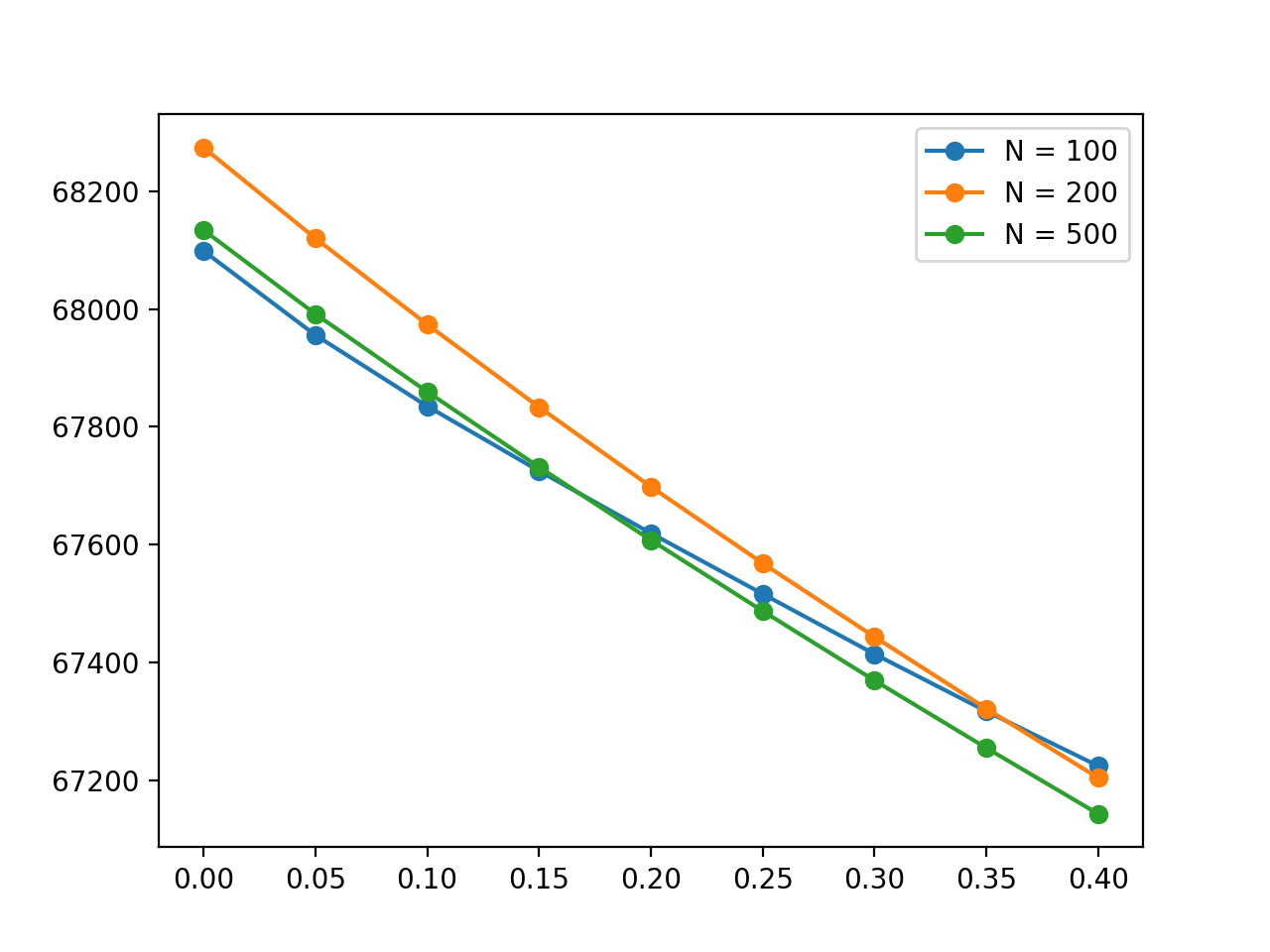} 
 \caption{$|\cS|=60$ and $B=20$.}
\end{subfigure}
\begin{subfigure}{.45\textwidth}
 \includegraphics[width=1\linewidth]{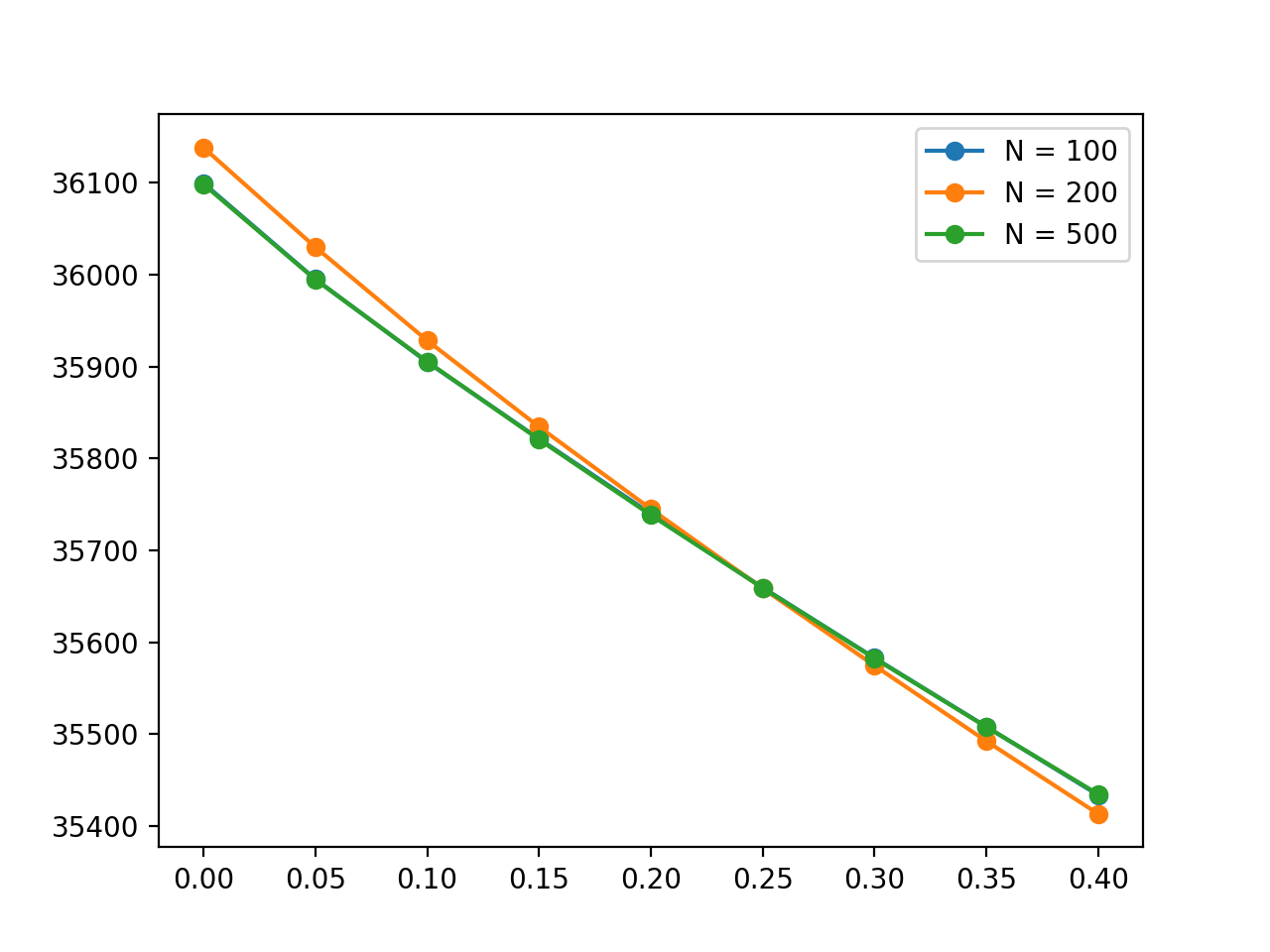} 
 \caption{$|\cS|=60$ and $B=10$}
\end{subfigure}
 \caption{The impact of total variation distance on the robust optimal objective of \eqref{opt:FL-DRO}.}
 \label{fig:obj-vs-TV}
\end{figure}

We further investigate the impact of the total variation distance $d$
on the robust optimal decision and the objective value. In this study,
we focus on the two instances $(|\cS|,B)=(60,20)$ and $(|\cS|,B)=(60,10)$. 
For a given sample size $N$ (with $N=100,200,500$), we tune the total variation distance 
$d$ from 0 to 0.4 with an increment 0.05 and investigate the change of optimal
objective of \eqref{opt:FL-DRO}. The results are shown in Figure~\ref{fig:obj-vs-TV}.
Obviously, the objective should decrease with $d$ as a larger $d$ value can allow
the worse scenarios get more weight and hence makes the objective decrease. 
With $d$ increasing from 0 to 0.4, the decrement of the objective is 1.5\% (1000 units)
for the case $B=20$, and 1.9\% (700 units) for the case $B=10$.
The relative decrement is about one magnitude smaller compared to the 10\% relative deviation in generating
samples for the demand. It is also observed that the optimal decision of \eqref{opt:FL-DRO}
does not change with the sample size or the total variation distance.

\begin{table}
\centering
{\scriptsize
\caption{Computational performance of solving 12 numerical instances of \eqref{opt:FL-DRO}.
In all instances, the total variation distance for the ambiguity set is set to be 0.2.
All instances are solved to optimality.}\label{tab:comp-perform}
\begin{tabular}{rr|rrr|rrr|rrr}
\hline\hline
 & & \multicolumn{3}{c|}{$N=100$} & \multicolumn{3}{c|}{$N=200$} & \multicolumn{3}{c}{$N=500$} \\
 \hline
$|\cS|$	&	$B$	&	Obj.($\times10^3$)	&	Iters. 	&	T(s)	&	Obj($\times10^3$)	&	Iters.	&	T(s)	&	Obj($\times10^3$)	&	Iters.	&	T(s)	\\
\hline
20	&	3	&	9.60	&	7	&	5.2	&	9.66	&	6	&	5.8	&	9.61	&	6	&	10.2	\\
20	&	5	&	15.57	&	9	&	2.3	&	15.57	&	7	&	4.8	&	15.6	&	8	&	11.2	\\
20	&	7	&	22.68	&	6	&	1.6	&	22.55	&	8	&	5.1	&	22.7	&	12	&	17.0	\\
\hline
40	&	5	&	17.52	&	10	&	8.2	&	17.55	&	11	&	17.8	&	17.53	&	10	&	43.4	\\
40	&	10	&	36.30	&	12	&	9.8	&	36.38	&	9	&	15.9	&	36.34	&	11	&	49.7	\\
40	&	15	&	54.44	&	11	&	9.2	&	54.46	&	11	&	19.1	&	54.38	&	12	&	55.7	\\
\hline
60	&	10	&	35.74	&	7	&	14.6	&	35.74	&	8	&	36.5	&	35.74	&	6	&	67.5	\\
60	&	15	&	54.64	&	8	&	17.4	&	54.70	&	10	&	46.8	&	54.70	&	10	&	120.2	\\
60	&	20	&	67.62	&	8	&	20.2	&	67.70	&	8	&	37.8	&	67.61	&	7	&	84.6	\\
\hline
80	&	15	&	60.62	&	13	&	61.3	&	60.47	&	7	&	74.3	&	60.56	&	6	&	158.1	\\
80	&	20	&	80.13	&	9	&	41.6	&	80.02	&	8	&	79	&	80.11	&	6	&	136	\\
80	&	25	&	98.3	&	14	&	70.6	&	98.29	&	14	&	127.5	&	98.38	&	10	&	220.5	\\
\hline
100	&	20	&	78.93	&	10	&	92.5	&	78.90	&	12	&	204.3	&	78.89	&	7	&	293.1	\\
100	&	25	&	100.69	&	13	&	120	&	100.54	&	11	&	185.5	&	100.53	&	15	&	637.6	\\
100	&	30	&	117.3	&	15	&	143.9	&	117.16	&	11	&	186.3	&	117.15	&	15	&	694.3	\\
\hline
200	&	40	&	173.24	&	12	&	1092.9	&	173.28	&	7	&	1184	&	173.41	&	11	&	4678.1	\\
200	&	50	&	215.74	&	13	&	1102.2	&	215.74	&	15	&	2549.9	&	215.78	&	9	&	3829.4	\\
200	&	60	&	256.8	&	14	&	1214.2	&	256.9	&	8	&	1360.1	&	256.78	&	13	&	5515.9	\\
\hline
\end{tabular}
}
\end{table}

\section{Concluding Remarks}
The distributionally-robust service center location problem investigated
in this research possesses decision-dependent demand induced naturally
from a maximum attraction principle and the number of opened service centers
in a neighborhood. The ambiguity set considered in this 
work is a decision-independent one defined using the total-variational distance
while many other metics such as the Wasserstein distance, Phi-divergences, etc.
can also be applied to define the ambiguity set with minor modification in the formulation. 
The modeling approach in this work decouples the endogenous impact (decision dependency) 
and distributional ambiguity. As a result, the model is highly computationally efficient
for mid- and large scale instances compared to models with an endogenous uncertainty 
(ambiguity) set. In practice, such decoupling is likely to be more data driven,
as there is a higher chance to identify data-based evidence that supports
decision-dependent demand than decision-dependent uncertainty of demand.

%
%
%
\begin{APPENDIX}{A cutting-plane algorithm for solving \eqref{opt:FL-DRO}}
\label{app:cut-plane-alg}
The model \eqref{opt:FL-DRO} is a special case of the distributionally-robust
two-stage stochastic program with a polyhedral ambiguity set.
\begin{algorithm}[H]
	{
	\caption{A cutting-plane algorithm for solving \eqref{opt:FL-DRO}.}
	\label{alg:cut-plane}
	\begin{algorithmic}[1]
	\State{Initialization: $\eta^*\gets\infty$, $y^*\gets 0^{|\F|}$, $n\gets 0$, $y^{(n)}\gets 0^{|\F|}$,
		 $\eta^{(n)}\gets\infty$ and $SolSet\gets\emptyset$.}
	\While{$y^{(n)}$ is not in $SolSet$.}
		\State{$SolSet\gets SolSet\cup\{y^{(n)}\}$.}
		\State{Set $n\gets n+1$.}
		\State{Solve the master problem \eqref{opt:master} for iteration $n$, and let $y^{(n)}$ be the optimal solution.}
		\State{Evaluate the value function $Q(y^{(n)},D^\omega)$ for each $\omega\in\Omega$.}
		\State{Solve the linear program \eqref{opt:LP-worst-prob} to get the worst-case probability measure $\mu^{(n)}$.}
		\State{Add the inequality \eqref{eqn:agg-ineq} to \eqref{opt:master}.}
	\EndWhile
	\State{Set $y^*\gets y^{(n)}$, and $\eta^*\gets\eta^{(n)}$.}
	\State{Return $y^*$ as the optimal solution and $\eta^*$ as the optimal objective.}
	\end{algorithmic}
	}
\end{algorithm}
\begin{theorem}[\citealt{bansal2018_decomp-alg-two-stage-dro}]
Algorithm~\ref{alg:cut-plane} terminates in a finite number of iterations and return an optimal solution
and optimal objective of \eqref{opt:FL-DRO}.
\end{theorem}
\end{APPENDIX}





\ACKNOWLEDGMENT{
The author is grateful to Dr. Liwei Zeng for a valuable discussion on getting an 
interpretation of the model. }

\bibliographystyle{informs2014} 
\bibliography{reference}



\end{document}